\documentclass[regno,12pt]{amsart}
\usepackage{amsmath,amsthm,amssymb,a4,graphics,color}
\usepackage[active]{srcltx}
\usepackage[usenames,dvipsnames]{pstricks}
\usepackage{epsfig}
\usepackage{pst-grad}
\usepackage{pst-plot} 
\usepackage{tikz}
\usepackage{hyperref}
\usepackage{dsfont}
\usepackage{float}

\newfam\Bbbfam 
\font\tenBbb=msbm10 
\font\sevenBbb=msbm7 
\font\fiveBbb=msbm5 
\textfont\Bbbfam=\tenBbb 
\scriptfont\Bbbfam=\sevenBbb 
\scriptscriptfont\Bbbfam=\fiveBbb

\newtheorem{theorem}{Theorem}[section] 
\newtheorem{lemma}[theorem]{Lemma} 
\newtheorem{prop}[theorem] {Proposition} 
\newtheorem{cor}[theorem]  {Corollary}

\theoremstyle{definition}

\newcommand{\dt}[1]{{\textrm{#1}}}

\def\1{{\mathchoice {1\mskip-4mu\mathrm l}      
{1\mskip-4mu\mathrm l} 
{1\mskip-4.5mu\mathrm l} {1\mskip-5mu\mathrm l}}}

\renewcommand{\qed}{\hfill\ensuremath{\square}}

\renewcommand{\d}{{\rm d}}

\newcommand{\e}   {{\operatorname e }}

\numberwithin{equation}{section}

\begin{document}
\title{A Spatial Model for Dormancy in Random Environment}
\author[Helia Shafigh]{}
\maketitle
\thispagestyle{empty}
\vspace{-0.5cm}

\centerline{\sc 
Helia Shafigh\footnote{WIAS Berlin, Mohrenstra{\ss}e 39, 10117 Berlin, Germany, {\tt shafigh@wias-berlin.de}}}
\renewcommand{\thefootnote}{}
\vspace{0.5cm}
\centerline{\textit{WIAS Berlin}}

\bigskip


\begin{abstract}
In this paper, we introduce a spatial model for dormancy in random environment via a two-type branching random walk in continuous-time, where individuals can switch between dormant and active states through spontaneous switching independent of the random environment. However, the branching mechanism is governed by a random environment which dictates the branching rates. We consider three specific choices for random environments composed of particles: (1) a Bernoulli field of immobile particles, (2) one moving particle, and (3) a Poisson field of moving particles. In each case, the particles of the random environment can either be interpreted as \emph{catalysts}, accelerating the branching mechanism, or as \emph{traps}, aiming to kill the individuals. The different between active and dormant individuals is defined in such a way that dormant individuals are protected from being trapped, but do not participate in migration or branching.

We quantify the influence of dormancy on the growth resp.\,survival of the population by identifying the large-time asymptotics of the expected population size. The starting point for our mathematical considerations and proofs is the parabolic Anderson model via the Feynman-Kac formula. Especially, the quantitative investigation of the role of dormancy is done by extending the Parabolic Anderson model to a two-type random walk.  
 \end{abstract}



\section{Introduction and main results}
\subsection{Biological Motivation} Dormancy is an evolutionary trait that has developed independently across various life forms and is particularly common in microbial communities. To give a definition, we follow \cite{blath} and refer to dormancy as the \emph{ability of individuals to enter a reversible state of minimal metabolic activity}. The collection of all dormant individuals within a population is also often called a \emph{seed-bank}. Maintaining a seed-bank leads to a decline in the reproduction rate, but it also reduces the need for resources, making dormancy a viable strategy during unfavourable periods. Initially studied in plants as a survival strategy (cf.\,\cite{cohen}), dormancy is now recognized as a prevalent trait in microbial communities with significant evolutionary, ecological, and pathogenic implications, serving as an efficient strategy to survive challenging environmental conditions, competitive pressure, or antibiotic treatment. However, it is at the same time a costly trait whose maintenance requires energy and
a sophisticated mechanisms for switching between active and dormant states. Moreover, the increased survival rate of dormant individuals must be weighed against their low reproductive activity. Despite its costs, dormancy still seems to provide advantages in variable environments. For a recent overview on biological dormancy and seed-banks we refer to \cite{dormancy}.

The existing stochastic models for dormancy can be roughly categorized into two approaches: population genetic models and population dynamic models. While the first approach assumes a constant population size and focusses on the genealogical implications of seed-banks, the latter is typically concerned with individual-based modelling through the theory of branching processes. Following a brief example in the book \cite{bookhaccou}, a two-type branching process (without migration) in a fluctuating random environment has been introduced in \cite{blath}, which served as a motivation for this paper. In \cite{blath}, the authors consider three different switching strategies between the two types (dormant and active), namely the \emph{stochastic (or: spontaneous; simultaneous) switching}, \emph{responsive switching} and \emph{anticipatory switching}. In the latter two strategies, individuals adapt to the fluctuating environment by selecting their state (dormant or active) based on environmental conditions via e.g.\,an increased reproduction activity during beneficial phases and a more extensive seed-bank during unfavourable ones (resp.\,vice versa). In contrast, the stochastic switching strategy, which remains unaffected by environmental changes, proves especially advantageous during catastrophic events, as it, with high probability, ensures the existence of dormant individuals, which may contribute to the survival of the whole population, when a severely adverse environment might eradicate all active ones. As an example, it is estimated that more than $80\%$ of soil bacteria are
metabolically inactive at any given time, forming extensive seed-banks of dormant individuals independent of the current conditions (cf. \cite{jl11} and \cite{ls18}). This makes the understanding of the stochastic switching strategy an interesting and important task. 

\subsection{Modelling Approach and Goals}  
The aim of this paper is to investigate the stochastic switching strategy in order to quantitatively compare the long-term behaviour of populations with and without this dormancy mechanism. 

Inspired by the Galton-Watson process with dormancy introduced in \cite{blath}, our first goal was to extend this model to a continuous-time spatial model on $\mathbb{Z}^d$. It is worth noting that spatial models for dormancy have already been considered in the setting of \emph{population genetics} (cf. \cite{spatialseedbanks}), where the population consists of different genetics types being inherited from parents to children. In such models, the total population size is fixed, so that the questions that arises are not about the extinction and survival of the whole population but rather about the evolution of the fraction of the different types. Notably, one of the goals in \cite{spatialseedbanks} is to determine criteria for \emph{co-existence} resp.\,\emph{clustering} of types in the limit of large population sizes. Another similar population genetics model, but this time with a (static) random environment, has been introduced in \cite{nandan}, in which the authors investigate the influence of dormancy again on co-existence and clustering. To the best of our knowledge, corresponding spatial models for dormancy in the setting of \emph{population size} models are still missing, such that the extension of the branching process in \cite{blath} seems to be a natural step. At the same time, there is a large repertoire of branching random walk models in random environment in the literature (cf.\,\cite{koenigsurvey} for a survey), even though none of them incorporates dormancy. Hence, by introducing a continuous-time spatial model with migration, branching resp.\,extinction driven by a random environment, as well as a Markovian switching between the two states \emph{active} and \emph{dormant}, we bridge the gap between  (non-spatial) two-type branching processes with dormancy on one side, and spatial branching random walks in random environments (without dormancy) on the other side of the existing literature. Especially, our main interest lies in a quantitative comparison of the population size of our model to those of existing branching random walk models without dormancy. 
\subsection{Description of the Model}
In our model, the population lives on $\mathbb{Z}^d$ and consists of two different types $i\in\{0,1\}$ of individuals, where we refer to $0$ as \emph{dormant} and to $1$ as \emph{active}. Let $\eta(x,i,t)$ be the number of individuals in spatial point $x\in\mathbb{Z}^d$ and state $i$ at time $t\geq 0$, which shall evolve in time according to the following rules:
\begin{itemize}
\item at time $t=0$, there is only one active individual in $0\in\mathbb{Z}^d$ and all other sites are vacant;
\item all individuals act independently of each other;
\item active individuals become dormant at rate $s_1\geq 0$ and dormant individuals become active at rate $s_0\geq 0$;
\item active individuals split into two at rate $\xi^+(x,t)\geq 0$ and die at rate $\xi^-(x,t)\geq 0$, depending on their spatial location $x$ and on time $t$, where both $\xi^+$ and $\xi^-$ are random fields;
\item active individuals jump to one of the neighbour sites with equal rate $\kappa\geq 0$;
\item dormant individuals do not participate in branching, dying or migration.
\begin{figure}[H]
\centering
\includegraphics[scale=0.8]{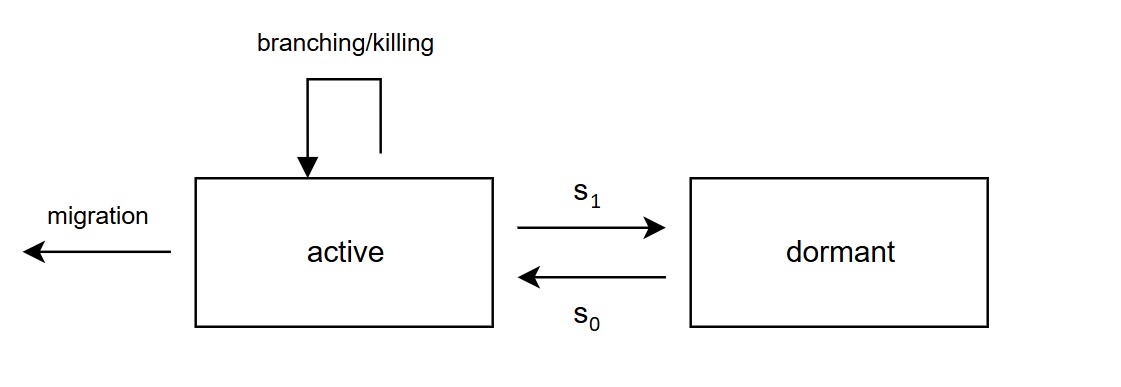}\caption{The evolution in every single point. Active individuals are subject to migration, branching and switching to dormant. Dormant individuals can only get active.}
\end{figure}
\end{itemize}
Write $\eta_0(x,i):=\eta(x,i,0)=\delta_{(0,1)}(x,i)$ and $\mathbb{P}_{\eta_0}$ for the corresponding probability measure with start in $\eta_0$. Then $(\eta,\mathbb{P}_{\eta_0})$ describes a Markov process on $\mathbb{N}^{\mathbb{Z}^d\times\{0,1\}}$. In the following, we abbreviate $\xi(x,t):=\xi^+(x,t)-\xi^-(x,t)$ for the \emph{balance} between branching and dying and refer to $\xi$ as the underlying \emph{random environment}. Let
\begin{align}\label{udef}
u(x,i,t):=\mathbb{E}_{\eta_0}[\eta(x,i,t)]
\end{align}
denote the expected number of individuals in $x\in\mathbb{Z}^d$ and state $i\in\{0,1\}$ at time $t$ with initial condition 
\begin{align*}
u(x,i,0)=\delta_{(0,1)}(x,i).
\end{align*}
Note, that the expectation is only taken over switching, branching and dying and not over the random environment $\xi$. If we average over $\xi$ as well, what we will denote in the following by $\left<\cdot\right>$, then we refer to
\begin{align*}
\left<u(x,i,t)\right>=\left<\mathbb{E}_{\eta_0}[\eta(x,i,t)]\right>
\end{align*}
as the \emph{annealed} number of individuals in $x\in\mathbb{Z}^d$ and in state $i\in\{0,1\}$ at time $t$. 

\subsection{Choices of the Random Environment}
In this paper we are going to consider random environments, which are built out of (interacting) particle systems. Regardless of the specific definition of each particle system, we will define the value of the random environment $\xi$ in point $x\in\mathbb{Z}^d$ at time $t$ as the number of particles present at this point, i.e.
\begin{align*}
\xi(x,t):=\# \text{ particles located at $x$ at time $t$}.
\end{align*}
In the following, we introduce three different particle systems, which we will consider as the underlying random environment throughout the paper:
\begin{enumerate}
\item \emph{Bournoulli field of immobile particles.} Place in each point $x\in\mathbb{Z}^d$ independently and with probability $p\in(0,1)$ one single particle which does not exhibit any movement. This results in a stationary field of immobile particles with a Bernoulli distribution over $\mathbb{Z}^d$, i.e.\,$\xi(x):=\xi(x,0)$ for all $x\in\mathbb{Z}^d$ and $t\geq 0$ and $\mathbb{P}(\xi(x)=1)=p$. 
\item \emph{One moving particle.} Here the random environment is dynamic and consists of one single particle starting in the origin and moving around according to a simple symmetric random walk $Y$ with total jump rate $2d\rho$. In other words,
\begin{align*}
\xi(x,t):=\delta_x(Y(t)), \qquad x\in\mathbb{Z}^d, t\geq 0. 
\end{align*}
\item \emph{Poisson field of moving particles.} At point $x\in\mathbb{Z}^d$, independently and according to a Poisson distribution with intensity $\nu$, place a random number of particles. The particles move independently of each other, each performing a simple symmetric random walk with total jump rate $2d\rho$. This setup generates a field of moving particles starting from a Poisson cloud. Then, we define the potential $\xi$ as
\begin{align*}
\xi(x,t):=\sum_{y\in\mathbb{Z}^d}\sum_{j=1}^{N_y}\delta_x(Y_j^y(t)), \qquad x\in\mathbb{Z}^d, t\geq 0, 
\end{align*}
where $N_y$ is a Poisson random variable with intensity $\nu>0$ for each $y\in\mathbb{Z}^d$, and $\{Y_j^y: y\in\mathbb{Z}^d, j=1,\cdots,N_y, Y_j^y(0)=y\}$ is the collection of independent random walks with total jump rate $2d\rho>0$. 
\end{enumerate}
Clearly, for each of the choices above, $\xi$ is a non-negative number, which results always in an positive balance between branching and killing. To allow for negative rates as well, we multiply $\xi$ with some factor $\gamma\in[-\infty,\infty)$ and will consider $\gamma\xi$ as the underlying random environment. Thus, each of our three choices can be either interpreted as a field of \emph{traps}, which corresponds to $\gamma<0$, or \emph{catalysts}, if $\gamma>0$. In the first case, active individuals will die with rate $|\gamma|$ if they encounter one of the traps, whereas they branch into two with rate $\gamma$ in the presence of catalysts in the latter case.   

\subsection{Results}
Recall the number of individuals $u(x,i,t)$ in point $x\in\mathbb{Z}^d$ and state $i\in\{0,1\}$ at time $t$, as defined in \eqref{udef}. The quantity we are interested in at most in the current paper is the \emph{annealed total number of individuals}
\begin{align}\label{grossu}
\left<U(t)\right>:=\sum_{x\in\mathbb{Z}^d}\sum_{i\in\{0,1\}}\left<u(x,i,t)\right>,
\end{align}
which turns into the \emph{annealed survival probability} up to time $t$ for $\gamma<0$. Our results concern the large-time asymptotics of $\left<U(t)\right>$ in case of both positive and negative $\gamma$ and for the three specific choices of the random environment mentioned above. 
Our first theorem quantifies the asymptotic behaviour of $\left<U(t)\right>$ if the environment consists of a Bernoulli field of particles:
\begin{theorem}\label{annasy1}
Let $\xi$ be chosen according to (1) and $d\geq 1$.
\begin{enumerate}
\item[(a)] If $\gamma=-\infty$, then the annealed survival probability converges to $0$, as $t\to\infty$, and obeys the asymptotics
\begin{align}\label{asy1-}
\log\left<U(t)\right>=-c_d\left|\log(1-p)\right|^{\frac{2}{d+2}}\left(\frac{\kappa 
s_0}{s_0+s_1}\right)^{\frac{d}{d+2}}t^{\frac{d}{d+2}}(1+o(1)), \quad t\to\infty,
\end{align}
for some constant $c_d$ depending only on the dimension $d$. 
\medskip
\item[(b)] If $\gamma\in(0,\infty)$, then the annealed number of individuals satisfies
\begin{align}\label{asy1+}
\lim_{t\to\infty}\frac{1}{t}\log\left<U(t)\right>=\gamma-s_1-\frac{(\gamma+s_0-s_1)^2-s_0s_1}{\sqrt{\gamma^2+2\gamma(s_0-s_1)+(s_0+s_1)^2}}.
\end{align}
\end{enumerate}
\end{theorem}
Our second theorem deals with the case of one moving particle:
\begin{theorem}\label{annasy2}
Let $\xi$ be chosen according to (2).
\begin{enumerate}
\item[(a)] If $\gamma\in(-\infty,0)$, then the annealed survival probability converges to $0$ in the dimensions $d\in\{1,2\}$, as $t\to\infty$, and satisfies the asymptotics
\medskip
\begin{align}\label{asy2-}
\left<U(t)\right>=\left\{\begin{array}{ll}\displaystyle\frac{2\sqrt{(s_0+s_1)(s_0(\rho+\kappa)+s_1\rho)}}{\sqrt{\pi}s_0|\gamma|}\frac{1}{\sqrt{t}}(1+o(1)), &d=1,\\[12pt]\displaystyle\frac{4\pi(s_1\rho+s_0(\rho+\kappa))}{s_0|\gamma|}\frac{1}{\log(t)}(1+o(1)), &d=2\end{array}\right.
\end{align}
as $t\to\infty$. 

In dimensions $d\geq 3$ the annealed survival probability admits the limit
In dimensions $d\geq 3$ the annealed survival probability admits the limit
\begin{align}\label{limitasy2-}
\lim_{t\to\infty}\left<U(t)\right>=1-\frac{|\gamma|G_d(0)}{\frac{s_0}{s_0+s_1}\left(\rho+\frac{s_0}{s_0+s_1}\kappa\right)+|\gamma|G_d(0)} \in(0,1),
\end{align}
where $G_d$ denotes the Green's function of a simple symmetric random walk with total jump rate $2d$. 
\\\item[(b)] If $\gamma\in(0,\infty)$, then for all $d\geq 1$, the annealed number of individuals satisfies
\begin{align}\label{asy2+}
\lim_{t\to\infty}\frac{1}{t}\log\left<U(t)\right>=\sup_{f\in\ell^2(\mathbb{Z}^d\times\{0,1\}),\|f\|_2=1}\left(A_1(f)-A_2(f)-A_3(f)\right)+\sqrt{s_0s_1}
\end{align}
where
\begin{align*}
A_1(f):=&\gamma f(0,1)^2,
\\A_2(f):=&\frac{1}{2}\sum_{i\in\{0,1\}}\sum_{x,y\in\mathbb{Z}^d, x\sim y}(i\kappa+\rho)(f(x,i)-f(y,i))^2,
\\A_3(f):=&\sum_{x\in\mathbb{Z}^d}\sqrt{s_0s_1}(f(x,1)-f(x,0))^2+\sum_{i\in\{0,1\}}\sum_{x\in\mathbb{Z}^d}s_i f(x,i)^2.
\end{align*}
\end{enumerate}
\end{theorem}
Finally, we establish the asymptotics of $\left<U(t)\right>$ for the third choice of the environment, namely a Poisson field of moving particles:
\begin{theorem}\label{annasy3}
Let $\xi$ be chosen according to (3).
\begin{enumerate}
\item[(a)] If $\gamma\in[-\infty,0)$, then the annealed survival probability converges to $0$ as $t\to\infty$ in all dimensions $d\geq 1$ and obeys the asymptotics
\medskip
\begin{align}\label{asy3-}
\log\left<U(t)\right>=\left\{\begin{array}{ll}\displaystyle-4\nu\sqrt{\frac{\rho s_0}{(s_0+s_1)\pi}}\sqrt{t}(1+o(1)), &d=1,\\[13pt]\displaystyle-4\nu\frac{\rho\pi s_0}{s_0+s_1}\frac{t}{\log\left(t\right)}(1+o(1)), &d=2,\\[13pt]\displaystyle-\lambda_{d,\gamma, \rho,\nu,s_0,s_1} t(1+o(1)), &d\geq 3,\end{array}\right.
\end{align}
as $t\to\infty$, for some constant $\lambda_{d,\gamma, \rho,\nu,s_0,s_1}>0$ depending on all the parameters.
\\\item[(b)] If $\gamma\in(0,\infty)$, then for all dimensions $d\geq 1$ the annealed number of individuals grows with double-exponential rate given by
\begin{align}\label{asy3+}
\lim_{t\to\infty}\frac{1}{t}\log\log\left<U(t)\right>=\sup_{f\in\ell^2(\mathbb{Z}^d),\|f\|_2=1}\left(\gamma f(0)^2-\frac{1}{2}\sum_{x,y\in\mathbb{Z}^d, x\sim y}\rho(f(x)-f(y))^2\right).
\end{align}
\end{enumerate}
\end{theorem}
\subsection{Relation to the Parabolic Anderson Model}
Recall the number of individuals $u(x,i,t)$ in point $x\in\mathbb{Z}^d$ and state $i\in\{0,1\}$ at time $t$ as defined in \eqref{udef}. It is already known (cf.\,\cite{baran}) that $u:\mathbb{Z}^d\times\{0,1\}\times[0,\infty)\to\mathbb{R}$ solves the partial differential equation
\medskip
\begin{align}\label{pamswitching}
\left\{\begin{array}{lllr}\frac{\d}{\d t}u(x,i,t)&=&i\kappa\Delta u(x,i,t) + Q u(x,i,t)+i\gamma\xi(x,t)u(x,i,t), &t>0, \\[12pt]u(x,i,0)&=&\delta_{(0,1)}(x,i),
\end{array}\right.
\end{align}
where
\begin{align*}
Q u(x,i,t):= s_i(u(x,1-i,t)-u(x,i,t))
\end{align*}
and $\Delta$ is the discrete Laplacian
\begin{align*}
\Delta f(x):=\sum_{y\in\mathbb{Z}^d, x\sim y}[f(y)-f(x)]
\end{align*}
acting on functions $f:\mathbb{Z}^d\to\mathbb{R}$, such that
\begin{align*}
\Delta u(x,i,t):=\sum_{y\in\mathbb{Z}^d, x\sim y}[u(y,i,t)-u(x,i,t)].
\end{align*}
We call \eqref{pamswitching} the \emph{parabolic Anderson model with switching}. If we consider a one-type branching random walk with only active individual evolving under the same evolution rules except of the switching mechanism, starting from one single particle in the origin, then it is well-known and has been shown in \cite{garmol} that the corresponding expected number of individuals solves the \emph{Parabolic Anderson model} (without switching)
\begin{align*}
\left\{\begin{array}{lllr}\frac{\d}{\d t}u(x,t)&=&\kappa\Delta u(x,t) + \gamma\xi(x,t)u(x,t), &t>0, x\in\mathbb{Z}^d \\[12pt]u(x,0)&=&\delta_0(x), &x\in\mathbb{Z}^d.
\end{array}\right.
\end{align*}
The parabolic Anderson model has been studied intensely during the past years and a comprehensive overview of results can be found in \cite{PAM}. One of the most powerful tools and often the starting point of the analysis of the PAM is the \emph{Feynman-Kac formula}
\begin{align}\label{fkwithout}
u(x,t)=\mathbb{E}_{x}^{X}\left[\exp\left(\int_0^t\gamma\xi(X(s),t-s)\,\d s\right)\delta_0(X(t))\right],
\end{align}
where $\mathbb{E}_{x}^{X}$ denotes the expectation over a simple symmetric random walk $X$ with start in $x$ and generator $\kappa\Delta$. In other words, the Feynman-Kac formula asserts that the time evolution of all individuals can be expressed as an expectation over one single individual moving around according to the same migration kernel and with a \emph{varying mass}, representing the population size. As we can see on the right hand-side of \eqref{fkwithout}, the mass of $X$ changes exponentially depending on the random environment $\xi$. Note, that if $\gamma<0$, then the right hand-side of \eqref{fkwithout} lies in $[0,1]$ and represents the \emph{survival probability} of a single individual up to time $t$. Now, since the Feynman-Kac formula is a powerful tool for the study of the parabolic Anderson model, it is only natural to pursue an analogous formulation in case of our two-type process with switching. To this end, let $\alpha=(\alpha(t))_{t\geq 0}$ be a continuous-time Markov process with state space $\{0,1\}$ and generator
\begin{align}\label{Q}
Qf(i):=s_i(f(1-i)-f(i))
\end{align}
for $f:\{0,1\}\to\mathbb{R}$. Conditioned on the evolution of $\alpha$, we define a continuous-time random walk $X=(X(t))_{t\geq 0}$ on $\mathbb{Z}^d$ which is supposed to stay still at a time $t$, if $\alpha(t)=0$, or perform a simple symmetric walk with jump rate $2d\kappa$, if $\alpha(t)=1$. In other words, the joint process $(X,\alpha)$ is the Markov process with the generator
\begin{align}\label{Lxa}
\mathcal{L}f(x,i):=i\kappa\sum_{y\sim x}(f(y,i)-f(x,i))+s_i(f(x,1-i)-f(x,i))
\end{align} 
for $x\in\mathbb{Z}^d$, $i,j\in\{0,1\}$ and a test function $f:\mathbb{Z}^d\times \{0,1\}\to\mathbb{R}$. Note, that the random walk $X$ itself is not Markovian due to the dependence on $\alpha$. Then, we call $(X,\alpha)$ a \emph{regime-switching random walk} (cf.\,\cite{switching} for the continuous-space version) and interpret $X$ as an individuals which is \emph{active} at time $t$, if $\alpha(t)=1$, and \emph{dormant} otherwise. Then, given a fixed realization of $\xi$, the formal solution of \eqref{pamswitching} is given by the \emph{Feynman-Kac formula}
\begin{align}\label{fkformula}
u(x,i,t)=\mathbb{E}_{(x,i)}^{(X,\alpha)}\left[\exp\left(\int_0^t\gamma\alpha(s)\xi(X(s),t-s)\,\d s\right)\delta_{(0,1)}(X(t),\alpha(t))\right],
\end{align}
where $\mathbb{E}_{(x,i)}^{(X,\alpha)}$ denotes the expectation over the joint process $(X,\alpha)$ starting in $(x,i)$ (cf.\,\cite{baran}). Thus, the study of our two-type branching process can be reduced to the analysis of only one individual with the same migration, branching and switching rates.

\subsection{Related Results}
The parabolic Anderson model without switching has been a topic of great interest during the past years. For a recent overview of results related to the classical model as well as many extensions we refer to \cite{PAM}. In this section we recall few results according to the parabolic Anderson model on $\mathbb{Z}^d$ which are most relevant for and most related to our models. Let us start with the case $\gamma<0$, i.e.\,with the case of a trapping random environment. The first model is best known as \emph{random walk among Bernoulli obstacles} and is the analogous version of our model (1) without switching, i.e.\,the time-independent potential $\xi$, which represents the \emph{traps} or \emph{obstacles}, is Bernoulli distributed with some parameter $p>0$. After getting trapped, the random walk with generator $\kappa\Delta$ dies immediately, which corresponds to the \emph{hard} trapping case $\gamma=-\infty$. It has been shown in \cite{antal} using a coarse graining technique known an \emph{method of
enlargement of obstacles} that the annealed survival probability up to time $t$ decays asymptotically as
\begin{align}\label{bernoulliwd}
\exp\left(-c_d(\log(1-p))^{\frac{2}{d+2}}(\kappa t)^{\frac{d}{d+2}}(1+o(1))\right), \quad t\to\infty,
\end{align}
with the same constant $c_d$ as in \eqref{asy1+}. In the setting of time-dependent potentials, the case of \emph{one moving trap} with generator $\rho\Delta$ with \emph{soft} killing ($\gamma\in(-\infty,0)$) has been studied in \cite{sw} for which it has been proven that the survival probability of a random walk with generator $\kappa\Delta$ among this mobile traps has the asymptotics
\begin{align}\label{schnitzler}
\left\{\begin{array}{ll}\displaystyle\frac{2\sqrt{\rho+\kappa}}{\sqrt{\pi}\gamma}\frac{1}{\sqrt{t}}(1+o(1)), &d=1,\\[12pt]\displaystyle\frac{4\pi(\rho+\kappa)}{|\gamma|\log(t)}(1+o(1)), &d=2,\\[12pt]\displaystyle 1-\frac{\gamma G_d(0)}{\rho+\kappa+|\gamma| G_d(0)}, &d\geq 3,\end{array}\right.
\end{align}
for $t\to\infty$ where $G_d$ denotes the Green's function of a random walk with generator $\Delta$. Hence, in dimension $d\in\{1,2\}$, the survival probability converges polynomially resp. logarithmically to zero, where the rate of convergence depends on all the parameters $\rho,\kappa$ and $\gamma$, whereas in dimensions $d\geq 3$ the survival probability converges to some number in $(0,1)$ and depends on the averaged total time spent in the origin, represented by the Green's function, as well. Finally, the case of a \emph{Poisson field of moving traps} with the random potential according to (3) has been investigated in \cite{drewitz}. Here, the survival probability of the random walk with jump rate $2d\kappa$ is known to decay asymptotically as
\begin{align}\label{resultsdrewitz}
\left\{\begin{array}{ll}\displaystyle\exp\left(-4\nu\sqrt{\frac{\rho}{\pi}}\sqrt{t}(1+o(1))\right), &d=1,\\[12pt]\displaystyle\exp\left(-4\nu\rho\pi\frac{t}{\log\left(t\right)}(1+o(1))\right), &d=2,\\[12pt]\displaystyle\exp(-\lambda_{d,\gamma,\rho,\nu} t(1+o(1))), &d\geq 3,\end{array}\right.
\end{align}
for some $\lambda_{d,\gamma,\rho,\nu}>0$, in both case of \emph{hard} and \emph{soft} trapping rates $\gamma\in[-\infty,0)$.\footnote{In \cite{drewitz} the authors have considered slightly different migration rates resulting in slightly different prefactors in the asymptotics, namely the normalized Laplacian $\frac{1}{2d}\Delta$ instead of $\Delta$. For better comparison, we stated \eqref{resultsdrewitz} in case of the non-normalized Laplacian $\Delta$.} That the survival probability, at least in dimensions $d\in\{1,2\}$, seems to be independent of the jump rate $\kappa$ of the individuals, is due to the fact that the asymptotics \eqref{resultsdrewitz} come from the behaviour that the individuals stay in the origin throughout the time, which corresponds to $\kappa=0$.

In the case $\gamma>0$ of catalysts, the analogous version without switching of our first model regarding Bernoulli distributed immobile catalysts is covered in \cite{bounded} as an example of a static bounded potential $\xi$. Although the results in \cite{bounded} are held very general, applying them to the Bernoulli potential yields that
\begin{align}\label{resultbiskup}
\lim_{t\to\infty}\frac{1}{t}\log\left<U(t)\right>=\gamma,
\end{align}
which come from the behaviour that the individuals find a region full of catalysts and stay there the whole time, such that they can continue branching throughout the time. 

For the case of one moving catalyst with jump rate $2d\rho$, is has been shown in \cite{movingcat} that the annealed number of individuals increases exponentially as well and the exponential growth rate it given by the variational formula
\begin{align}\label{resultheydenreich}
\lim_{t\to\infty}\frac{1}{t}\log\left<U(t)\right>=\sup_{f\in\ell^2(\mathbb{Z}^d),\|f\|_2=1}\left(\gamma f(0)^2-\frac{1}{2}\sum_{x,y\in\mathbb{Z}^d, x\sim y}(\kappa+\rho)(f(x)-f(y))^2\right).
\end{align} 
A similar result but with double-exponential growth has been proven in \cite{intermittency} for the case of a Poisson field of moving catalysts, where the rate
\begin{align}\label{resultgaertner}
\lim_{t\to\infty}\frac{1}{t}\log\log\left<U(t)\right>=\sup_{f\in\ell^2(\mathbb{Z}^d),\|f\|_2=1}\left(\gamma f(0)^2-\frac{1}{2}\sum_{x,y\in\mathbb{Z}^d, x\sim y}\rho(f(x)-f(y))^2\right).
\end{align}
has been shown to be finite in all dimensions $d\geq 1$. The absence of $\kappa$ is due to the fact that it is most favourable for the population growth if the individuals stay immobile, which corresponds to $\kappa=0$. Moreover, is has been shown in \cite{intermittency} that there are cases, where already the quantity $\frac{1}{t}\log\left<U(t)\right>$ converges to a finite limit, namely in the transient dimensions $d\geq 3$ under the assumption that $0\leq\frac{\gamma}{\rho}\leq G_d(0)^{-1}$, where $G_d(0)$ denotes the Green's function again. In these regimes \eqref{resultgaertner} therefore converges to $0$.
\subsection{Discussion}
In this section, we discuss the extent to which the stochastic dormancy strategy, as defined in our model, affects the long-time dynamics of the population. Let us start with the case of catalytic random environments. 
Recall from \eqref{resultbiskup} that if the environment is chosen according to a static Bernoulli field, then in the analogous model without dormancy the population grows exponentially fast in time $t$ with rate $\gamma$. Theorem 1.1.(b) asserts that the population growth still occurs exponentially in $t$ when dormancy is incorporated; however, the growth rate is no longer $\gamma$ any more but rather a smaller constant, as seen from \eqref{asy1+}. Indeed, an easy calculation shows that 
\begin{align*}
s_1+\frac{(\gamma+s_0-s_1)^2-s_0s_1}{\sqrt{\gamma^2+2\gamma(s_0-s_1)+(s_0+s_1)^2}}>0
\end{align*}
for all $\gamma,s_0,s_1>0$, such that the growth rate is strictly decreased after incorporating the stochastic dormancy strategy. It is worth mentioning that this effect comes from the probability for \emph{large deviations} of the time spent in the active state. Indeed, as we will see later in the proofs, \eqref{asy1+} can also be expressed as
\begin{align*}
\lim_{t\to\infty}\frac{1}{t}\log\left<U(t)\right>=\gamma-\inf_{[a\in[0,1]}\{I(a)+\gamma(1-a)\},
\end{align*}
where $I$ is a non-negative and strictly convex function defined as
\begin{align}\label{definitionI}
I(a)=(\sqrt{s_1a}-\sqrt{s_0(1-a)})^2.
\end{align}
Next, note that the population only grows in those time intervals in which the individuals are active. Now, on one hand, we have a \emph{law of large numbers} for the fraction of time spent in the active state, which asserts that the average time proportion each individual spends in the active state equal $s_0/(s_0+s_1)$ with probability one. On the other hand, as we will prove in Section 2, there is a \emph{large deviation principle} asserting that for any other fraction of time $a\in[0,1]$, the probability of spending $a\cdot t$ time units in the active state up to time $t$ decays exponentially in $t$ for $t\to\infty$, where the decay rate depends on the proportion $a$. This decays rate is nothing but the function $I$ defined in \eqref{definitionI}. Hence, \eqref{asy1+} tells us that the original growth rate $\gamma$ is decreased by $\gamma$ multiplied with the proportion $1-a$ of time spent in the dormant state on one hand, as there is no reproduction in this case, and by $I(a)$ on the other hand, as it represents the probabilistic cost to spend exactly the proportion $a$ of time in the active state. At the end, this probabilistic cost has to be weighed against the positive contribution to reproduction, such that $I(a)+\gamma(1-a)$ has to be optimized over $a$. 

Continuing with one moving catalyst and comparing \eqref{asy2+} to \eqref{resultheydenreich}, we see that the population again grows exponentially in $t$ and the rate is affected by all the involved mechanisms. While the first term in the variational formula \eqref{resultheydenreich} shows that branching with rate $\gamma$ occurs whenever the distance between individuals and catalytic particles is equal to $0$, the term $\gamma f(0,1)^2$ appearing in \eqref{asy3+} has the interpretation that in case of dormancy, individuals can only branch with rate $\gamma$ if, first, the distance between them and the catalyst is equal to $0$, and second, if they are in the active state $1$. To highlight another difference, we see that the only probabilistic cost that appears in the variational formula \eqref{resultheydenreich} is the one coming from the movement of the individuals (with rate $2d\kappa$) as well as of the catalyst (with rate $2d\rho$), whereas in \eqref{asy2+} besides the movements appearing  in the term $A_2$, also the exchange between states, represented by $A_3$, has to be taken into account. The additional term $\sqrt{s_0s_1}$ comes from a change of measure , which will be clarified in Section 2. 

In case of a Poisson field of moving catalysts, we see that our asymptotics \eqref{asy3+} is on a double-exponential scale and equals the variational formula \eqref{resultgaertner} for the corresponding model without dormancy. In other words, although the stochastic dormancy strategy slows down the population growth due to the lack of reproduction in dormant phases in the first two choices of the environment, this inactivity does not seem to influence the population growth at all, if the moving catalysts start from a Poisson cloud. As will be revealed in Section 5, this is due to the fact that if the individuals manage to find favourable regions with a high density of catalysts, then the reproduction rate in the active state is on such a high scale, namely double-exponentially in time $t$, that the exponential probabilistic cost to stay active is negligible in comparison to the high positive outcome. Thus, the variational formula \eqref{asy3+} does not take dormancy into account and depends only on the branching rate and movement of the catalysts. 
As we will see later in the proofs, the independence of $\kappa$ arises from the fact that it is most favourable for the population growth if the individuals stay immobile, which corresponds to $\kappa=0$, and matches the behaviour in the model without dormancy, as mentioned in the last section. 

Next, we discuss the case $\gamma<0$ of trapping environments. Our results concerning this case can be summarized more briefly, since they share a similarity regarding the dependence on dormancy: At least in dimensions, in which we have an explicit expression for the asymptotic survival probability, we see that the latter is increased after incorporating the stochastic dormancy strategy in comparison to the corresponding models without dormancy, and is monotone in the average time $s_1/(s_0+s_1)$ spent in the dormant state. Moreover, setting $s_1=0$ yields exactly the same asymptotics as for the models without dormancy, making our results a generalization for arbitrary $s_1\geq 0$. This is immediate clear for the Bernoulli field of immobile traps in all dimensions, as comparing \eqref{asy1-} and \eqref{bernoulliwd} shows. As will be addressed in the proof, formula \eqref{asy1-} also indicates that the higher survival probability results from a \emph{time-change}, since individuals can only move towards the immobile traps during active phases. Therefore, if we only take into account those time intervals in which the individuals are mobile, which in average accounts for a proportion of $s_0/(s_0+s_1)$ of the whole time due to the law of large numbers, then \eqref{bernoulliwd} translates into \eqref{asy1-}. Note that here the law of large numbers dictates the behaviour of the time spent in the active state and not the large deviation principle, since the scale $t^{d/(d+2)}$ is much smaller than the large deviation scale $t$. 

The law of large number seems to play a role also in the case of one moving catalyst. However, comparing \eqref{asy2-} and \eqref{schnitzler} demonstrates that in this case the positive effect of dormancy on the survival probability does not only come from a time-change. Rewriting the pre-factor
\begin{align*}
\frac{2\sqrt{(s_0+s_1)(s_0(\rho+\kappa)+s_1\rho)}}{\sqrt{\pi}s_0|\gamma|}=\frac{2\sqrt{\frac{s_0}{s_0+s_1}\kappa+\rho}}{\sqrt{\pi}\frac{s_0}{s_0+s_1}|\gamma|}
\end{align*}
of the polynomial asymptotics \eqref{asy2-} in dimension $d=1$ and comparing it to \eqref{schnitzler} suggests that, although the time-change is still present as a pre-factor of the jump rate $\kappa$, the killing rate $\gamma$ is reduced as well, since the individuals are again only a proportion of $s_0/(s_0+s_1)$ of the time active and therefore vulnerable to the traps. For $d=2$, we can see both effects as well. Also the monotonicity of the survival probability in the average dormancy proportion $s_1/(s_0+s_1)$ becomes evident  through a simple calculation. 
We see similar  effects also in case of a Poisson field of moving traps by comparing \eqref{asy3-} to \eqref{resultsdrewitz}, where the reduction of the exponential decay rate by factor $\sqrt{s_0/(s_0+s_1)}$ in dimension $1$ and by factor $s_0/(s_0+s_1)$ in dimension $2$ comes from the law of large numbers as well, whereas in the dimensions $d\geq 3$ the large deviation principle dictates the asymptotics due to the joint time scale $t$. The (surprising) independence of the survival probability of the killing rate $\gamma$ as well as the jump rate $\kappa$ has been discussed in \cite{drewitz} for the corresponding model without dormancy and underlies the same reasons in our case. 

To summarize, for $\gamma>0$ the number of individuals up to time $t$ is either unchanged or reduced due to dormancy in our models, while for $\gamma<0$, at least in those cases in which an explicit expression is given, the survival probability is increased with dormancy and is monotone in the average time spent in the dormant state. 
\begin{figure}[H]
\centering
\includegraphics[scale=0.8]{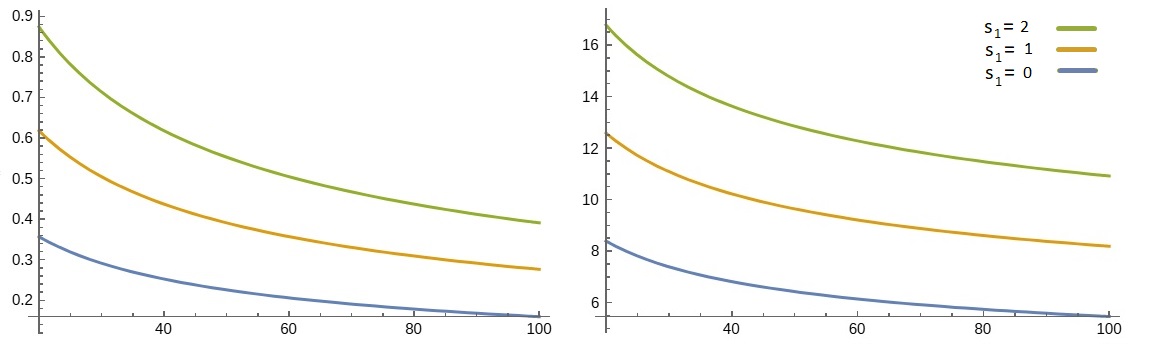}\caption{Asymptotics of the annealed survival probability in case of one moving trap in dimension $1$ (left) and $2$ (right) as a function of $t$ for $\rho=\kappa=s_0=\gamma=1$ and different choices of $s_1$.}
\end{figure}
\begin{figure}[H]
\centering
\includegraphics[scale=0.8]{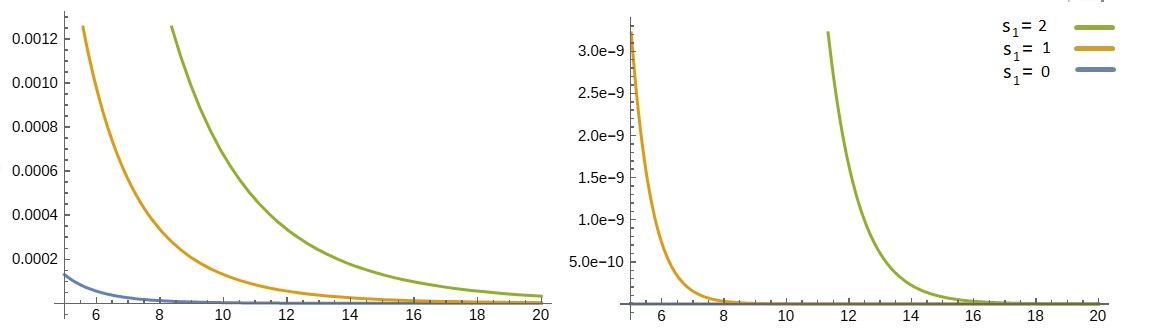}\caption{Asymptotics of the annealed survival probability in case of a Poisson field of moving traps in dimension $1$ (left) and $2$ (right) as a function of $t$ for $\rho=\kappa=s_0=\gamma=1$ and different choices of $s_1$.}
\end{figure}
\subsection{Outline}
The rest of this paper is organized as follows. In section 2 we will convert our branching model into a switching random walk via the Feynman-Kac formula in order to obtain a more convenient representation of $\left<U(t)\right>$. Further, we will collect some results related to the distribution and large deviations of the local times of $\alpha$. Section 3, 4 and 5 respectively are dedicated to the proofs of our main Theorems 1.1, 1.2 and 1.3 respectively.
\section{Preparatory facts}
\subsection{Feynman-Kac formula}
As discussed in the introduction, our dormancy model is motivated by population dynamics and initially defined as a two-type branching random walk with Markovian switching between the types. However, all our proofs and considerations are based on the Feynman-Kac formula \eqref{fkformula}, which serves as the cornerstone for the subsequent steps throughout the remainder of this paper. Note that our choices of our dynamic random environments (2) and (3) are reversible in time, in the sense that $(\xi(\cdot,t))_{0\leq t\leq T}$ is equally distributed to $(\xi(\cdot,T-t))_{0\leq t\leq T}$, for all $T>0$. Hence, taking the expectation with respect to $\xi$ and changing the order of integrals, we can write
\begin{align*}
\left<U(t)\right>&=\sum_{x\in\mathbb{Z}^d}\sum_{i\in\{0,1\}}\mathbb{E}_{(x,i)}^{(X,\alpha)}\left[\left<\exp\left(\int_0^t\gamma\alpha(s)\xi(X(s),t-s)\,\d s\right)\delta_{(0,1)}(X(t),\alpha(t))\right>\right]
\\&=\sum_{x\in\mathbb{Z}^d}\sum_{i\in\{0,1\}}\mathbb{E}_{(x,i)}^{(X,\alpha)}\left[\left<\exp\left(\int_0^t\gamma\alpha(s)\xi(X(s),s)\,\d s\right)\delta_{(0,1)}(X(t),\alpha(t))\right>\right].
\end{align*}
Moreover, we can change the starting and end point of the path $(X,\alpha)$ to obtain
\begin{align}\label{repu}
\left<U(t)\right>&=\sum_{x\in\mathbb{Z}^d}\sum_{i\in\{0,1\}}\mathbb{E}_{(0,1)}^{(X,\alpha)}\left[\left<\exp\left(\int_0^t\gamma\alpha(s)\xi(X(s),s)\,\d s\right)\delta_{(x,i)}(X(t),\alpha(t))\right>\right]\nonumber
\\&=\mathbb{E}_{(0,1)}^{(X,\alpha)}\left[\left<\exp\left(\int_0^t\gamma\alpha(s)\xi(X(s),s)\,\d s\right)\right>\right]\nonumber
\\&=\left<\mathbb{E}_{(0,1)}^{(X,\alpha)}\left[\exp\left(\int_0^t\gamma\alpha(s)\xi(X(s),s)\,\d s\right)\right]\right>.
\end{align}
Especially, $U(t)$ can also be interpreted as the solution $\tilde{u}$ of the partial differential equation
\begin{align*}
\left\{\begin{array}{lllr}\frac{\d}{\d t}\tilde{u}(x,i,t)&=&i\kappa\Delta \tilde{u}(x,i,t) + s_i (\tilde{u}(x,1-i,t)-\tilde{u}(x,i,t))+i\gamma\xi(x,t)\tilde{u}(x,i,t), &t\geq 0, \\[12pt]\tilde{u}(x,i,0)&=&1
\end{array}\right.
\end{align*}
in point $(0,1)$, which differs from \eqref{pamswitching} only in the initial distribution. From now on, we will work with the representation \eqref{repu}. Note, that \eqref{repu} also holds for the static choice (1), as $\xi$ does not depend on time. 
\subsection{Large deviation principle for the relative active times}
Let
\begin{align*}
L_t(i):=\int_0^t\delta_i(\alpha(s))\d s
\end{align*}
denote the local times of $\alpha$ in state $i\in\{0,1\}$ up to time $t$, i.e. the total time till $t$ spent in state $i$. The goal of this section is to prove a large deviation principle for the normalized local times $(\small\frac{1}{t}L_t(i))_{t\geq 0}$, i.e. for the proportion of time an individual spends in the state $i$. Where such a large deviation principle for normalized local times is already well-known in the literature for discrete-space Markov processes with \emph{symmetric} transition rates (cf.~\cite[Theorem 3.6.1 and Remark 3.6.4]{koenigldp}), the corresponding principle in case of \emph{asymmetric} rates is still missing. In our case, we can obtain the probability for large deviations directly by computing the exact distribution of the local times, which will be done in the following lemma:
\begin{lemma}\label{distlocaltimes}
For all $s_0,s_1>0$ the probability density function of the local times $(L_t(1))_{t>0}$ of $\alpha$ in state $1$ is given by
\begin{align*}
\mathbb{P}(L_t(1)\in\d y)=&s_1\e^{-s_0t-(s_1-s_0)y}\left(\sum_{k=0}^{\infty}\frac{(s_0s_1y(t-y))^k}{k!k!}\left(\frac{s_0y}{k+1}+1\right)\right).
\end{align*}
\end{lemma}
\begin{proof}[Proof]
Let $N(t)$ be the number of jumps of the Markov chain $\alpha$ up to time $t>0$ when $\alpha$ starts in state $1$. We denote by $a_i$, $i\in\mathbb{N}$, the waiting times of transitions from $0$ to $1$ and by $b_i$ those from $1$ to $0$ such that $a_i$ resp.~ $b_i$ are independent and exponentially distributed with parameter $s_0$ resp.~ $s_1$. If $N(t)$ is even, $\alpha$ will be in state $1$ again after the last jump before $t$. Then
\begin{align*}
\mathbb{P}(L_t(0)&\in\d y, N(t)=2k)
\\&=\mathbb{P}\left(\sum_{i=1}^ka_i\in\d y, \sum_{i=1}^ka_i+\sum_{i=1}^kb_i<t, b_{k+1}>t-\left(\sum_{i=1}^ka_i+b_i\right)\right)
\\&=\int_{x=0}^{t-y}\mathbb{P}\left(\sum_{i=1}^kb_i\in\d x, \sum_{i=1}^ka_i\in\d y, b_{k+1}>t-x-y\right)
\\&=\int_0^{t-y}\frac{s_1^kx^{k-1}e^{-s_1x}}{(k-1)!}\frac{s_0^ky^{k-1}\e^{-s_0y}}{(k-1)!}\e^{-s_1(t-x-y)}\,\d x\,\d y
\\&=\frac{(s_0s_1)^k}{(k-1)!k!}\e^{-s_1t-(s_0-s_1)y}(t-y)^ky^{k-1}
\end{align*}
for $y\in[0,t]$ and therefore
\begin{align*}
\mathbb{P}(L_t(1)\in\d y, N(t)=2k)=\frac{(s_0s_1)^k}{(k-1)!k!}\e^{-s_0t-(s_1-s_0)y}(t-y)^{k-1}y^{k}.
\end{align*}
In case of an odd number of jumps, where $\alpha$ is in state $0$ after the last jump before $t$, the joint distribution of $L_t(1)$ and $N(t)$ reads
\begin{align*}
\mathbb{P}(L_t(1)\in\d y, N(t)=2k+1)=&\int_{x=0}^{t-y}\mathbb{P}\left(\sum_{i=1}^{k+1}b_i\in\d y, \sum_{i=1}^ka_i\in\d x, a_{k+1}>t-x-y\right)
\\=&\frac{s_0^ks_1^{k+1}}{k!k!}\e^{-s_0t-(s_1-s_0)y}(t-y)^ky^k.
\end{align*}
The claim follows after summing over all $k\in\mathbb{N}$. 
\end{proof}
Making use of the exact distribution of the local times of $\alpha$, we are able to establish the following large deviation principle:
\begin{cor}[LDP for local times of $\alpha$]\label{LDP}
For all choices of the transition rates $s_0,s_1>0$, the normalized local times $\left(\frac{1}{t}L_t(1)\right)_{t\geq 0}$ of $\alpha$ in state $1$ satisfy a large deviation principle on $[0,1]$ with rate function $I:[0,1]\to\mathbb{R}$ given by
\begin{align}\label{I}
I(a)=(\sqrt{s_1a}-\sqrt{s_0(1-a)})^2.
\end{align}
\end{cor}
\begin{proof}
This follows immediately from Lemma \ref{distlocaltimes}. Indeed, observe that for $x,y>0$,
\begin{align*}
\sum_{k=0}^{\infty}\frac{x^k}{(k!)^2}\left(\frac{y}{k+1}+1\right)&=\sum_{k=0}^{\infty}\frac{x^k}{(k!)^2}+y\sum_{k=0}^{\infty}\frac{x^k}{(k!)^2(k+1)}
\\&=I_0(2\sqrt{x})+y\int_0^1I_0(2\sqrt{x u})\d u,
\end{align*}
where 
$I_0(z)=\sum_{k=0}^{\infty}\frac{\left(\frac{1}{4}z^2\right)^k}{k!k!}$ is the modified Bessel function with parameter $0$, which can also be written as
\begin{align}\label{besselrep}
I_0(z)=\frac{1}{\pi}\int_0^{\pi}\e^{z\cos(\theta)}\,\d\theta.
\end{align}
Hence, for all $t>0$ and $a\in[0,1]$,
\begin{align*}
\mathbb{P}(L_t(1)&\in\d (at))=\e^{-s_0t-(s_1-s_0)at}
\\&\times\left(\frac{1}{\pi}\int_0^\pi\e^{2t\sqrt{s_0s_1a(1-a)}\cos(\theta)}\d\theta+\frac{s_0at}{\pi}\int_0^1\int_0^\pi\e^{2t\sqrt{s_0s_1a(1-a)u}\cos(\theta)}\d\theta\d u\right).
\end{align*}
Applying the method of Laplace for integrals (cf.~ \cite[Corollary 1.3.2]{koenigldp}) yields
\begin{align*}
\lim_{t\to\infty}\frac{1}{t}\log\left(\int_0^{\pi}\e^{t2\sqrt{s_0s_1a(1-a)}\cos(\theta)}\,\d\theta\right)&=\max_{\theta\in[0,\pi]}2\sqrt{s_0s_1a(1-a)}\cos(\theta)\\&=2\sqrt{s_0s_1a(1-a)}
\end{align*}
and analogously
\begin{align*}
\lim_{t\to\infty}\frac{1}{t}\log\left(\frac{s_0at}{\pi}\int_0^1\int_0^\pi\e^{2t\sqrt{s_0s_1a(1-a)u}\cos(\theta)}\d\theta\d u\right)&=\max_{\substack{u\in[0,1],\\\theta\in[0,\pi]}}2\sqrt{s_0s_1a(1-a)u}\cos(\theta)
\\&=2\sqrt{s_0s_1a(1-a)}.
\end{align*}
Consequently,
\begin{align*}
I(a):=&-\lim_{t\to\infty}\frac{1}{t}\log \mathbb{P}^{\alpha}_1(L_t(1)\in\d (at))
\\=&s_1a+s_0(1-a)-2\sqrt{s_0s_1a(1-a)}=(\sqrt{s_1a}-\sqrt{s_0(1-a)})^2.
\end{align*}
\end{proof}
Note, that the rate function $I$ is strictly convex and, as an easy calculations shows, attends its unique minimizer at $a_{\text{min}}=\frac{s_0
}{s_0+s_1}$. This implies a law of large numbers for the proportion of time spent in the active state, i.\,e.\,
\begin{align*}
\lim_{t\to\infty}\frac{1}{t}L_t(1)=\frac{s_0}{s_0+s_1} 
\end{align*}
almost surely. Further, if we choose $s:=s_0=s_1$ equally, the rate function becomes
\begin{align*}
I(a)=s(\sqrt{a}-\sqrt{1-a})^2
\end{align*}
which is the well-known large deviation rate function in the case of symmetric transition rates (cf.\,\cite[Theorem 3.6.1 and Remark 3.6.4]{koenigldp}).
The next lemma shows that if we look at the exponential moments of the normalized local times $\left(\frac{1}{t}L_t(1)\right)_{t\geq 0}$ of $\alpha$ with a \emph{smaller} exponential rate than $t$, then the best $\left(\frac{1}{t}L_t(1)\right)_{t\geq 0}$ can do in order to maximize the exponent, is to take its long-term average:
\begin{lemma}\label{varadhanversion}
Let $(f(t))_{t\geq 0}$ be a sequence of positive real numbers with $\lim_{t\to\infty}f(t)=\infty$ and $\lim_{t\to\infty}\frac{f(t)}{t}=0$. Then, for any continuous and bounded function $F\colon[0,1]\to\mathbb{R}$,
\begin{align*}
\lim_{t\to\infty}\frac{1}{f(t)}\log\mathbb{E}_{1}^{\alpha}\left[\e^{f(t)F\left(\frac{1}{t}L_t(1)\right)}\right]=F\left(\frac{s_0}{s_0+s_1}\right).
\end{align*}
\end{lemma}
\begin{proof}
For the lower bound, fix $\delta>0$ and let $G$ be some open ball around $x_0:=\frac{s_0}{s_0+s_1}$ such that $F(x_0)-\delta\leq\inf_G F(x_0)\leq  F(x_0)+\delta$. Then we use the law of large numbers for the sequence $\left(\frac{1}{t}L_t(1)\right)_{t\geq 0}$ to obtain
\begin{align*}
\liminf_{t\to\infty}\frac{1}{f(t)}\log\mathbb{E}_{1}^{\alpha}\left[\e^{f(t)F\left(\frac{1}{t}L_t(1)\right)}\right]\geq&\liminf_{t\to\infty}\frac{1}{f(t)}\log\left(\e^{f(t)\inf_{a\in G}F(a)}\mathbb{P}_{1}^{\alpha}\left(\frac{1}{t}L_t(1)\in G\right)\right)
\\\geq& \inf_{a\in G}F(a)+\liminf_{t\to\infty}\frac{1}{f(t)}\log\mathbb{P}_{1}^{\alpha}\left(\frac{1}{t}L_t(1)\in G\right)
\\\geq& F(x_0)-\delta.
\end{align*}
For the upper bound, the boundedness of $F$ by some $M>0$ gives
\begin{align*}
\mathbb{E}_1^{\alpha}\left[\e^{f(t)F\left(\frac{1}{t}L_t(1)\right)}\right]\leq& \e^{f(t)(F(x_0)+\delta)}\mathbb{P}_1^{\alpha}\left(\frac{1}{t}L_t(1)\in\bar{G}\right)+\e^{f(t)M}\mathbb{P}_1^{\alpha}\left(\frac{1}{t}L_t(1)\in G^c\right).
\end{align*}
Now, as
\begin{align*}
\limsup_{t\to\infty}\frac{1}{f(t)}\log\mathbb{P}_1^{\alpha}\left(\frac{1}{t}L_t(1)\in G^c\right)=-\limsup_{t\to\infty}\frac{t}{f(t)}\inf_{a\in G^c}I(a)=-\infty,
\end{align*}
where $I$ is the rate function defined in \eqref{I}, and using the law of large numbers again as well as the method of Laplace for sums, we deduce that
\begin{align*}
\limsup_{t\to\infty}\frac{1}{f(t)}\log\mathbb{E}_1^{\alpha}\left[\e^{f(t)F\left(\frac{1}{t}L_t(1)\right)}\right]\leq \max\left\{F(x_0)+\delta, -\infty\right\}=F(x_0)+\delta.
\end{align*}
The claim follows after letting $\delta\to 0$. 
\end{proof}
\subsection*{Change of measure for $\alpha$}
One of the proof methods we will use in section 4 to obtain the representation \eqref{asy2+} is the Perron-Frobenius spectral theory for bounded self-adjoint operators, which we would like to apply to the generator
\begin{align}\label{x-y}
\bar{\mathcal{L}}f(z,i)=\sum_{y\sim z}(i\kappa+\rho)(f(y,i)-f(z,i))+s_i(f(z,1-i)-f(z,i))
\end{align} 
of a Markov process $(Z,\alpha)$, where $f\colon\mathbb{Z}^d\times\{0,1\}\to\mathbb{R}$ is a suitable test function. However, this may be a problem, as the matrix $Q$ defined in \eqref{Q} is not symmetric and therefore the generator \eqref{x-y} of $(Z,\alpha)$ is not self-adjoint. In order to fix this, we will use a result from \cite{girsanovmarkov} which we formulate here for the convenience of the reader:
\begin{lemma}\label{girsanov}
Let $\alpha=(\alpha(t))_{t\geq 0}$ be any Markov process on a finite state space $\mathcal{M}$ with transition rates $q_{ij}$ from state $i\in\mathcal{M}$ to $j\in\mathcal{M}$. For a positive function $h\colon\mathcal{M}\to(0,\infty)$, let $\tilde{\alpha}=(\tilde{\alpha}(t))_{t\geq 0}$ be another Markov process on $\mathcal{M}$ defined on the same filtered space $(\Omega, (\mathcal{F}_t)_{t\geq 0}))$ with transition rates $\tilde{q}_{ij}$ given by
\begin{align*}
\tilde{q}_{ij}=q_{ij}\frac{h(j)}{h(i)}
\end{align*}
for $i\neq j$ and $\tilde{q}_{ii}=-\sum_{k\neq i}q_{ik}\frac{h(k)}{h(i)}$. Denote by $\mathbb{P}^{\alpha}$ resp.\,$\mathbb{P}^{\tilde{\alpha}}$ the distribution of $\alpha$ resp.\,$\tilde{\alpha}$. Then $\mathbb{P}^{\alpha}$ is absolutely continuous with respect to $\mathbb{P}^{\tilde{\alpha}}$ with the Radon-Nikodym derivative 
\begin{align*}
\frac{\d\mathbb{P}^{\alpha}}{\d\mathbb{P}^{\tilde{\alpha}}}|_{\mathcal{F}_t}=\frac{h(\tilde{\alpha}(t))}{h(\tilde{\alpha}(0))}\exp\left(-\int_0^t\frac{\tilde{Q}h(\tilde{\alpha}(s))}{h(\tilde{\alpha}(s))}\,\d s\right),
\end{align*}
where the operator $\tilde{Q}$ is defined as
\begin{align*}
\tilde{Q}f(i)=\tilde{q}_{ij}(f(j)-f(i))
\end{align*}
for $f:\mathcal{M}\to(0,\infty)$. 
\end{lemma}
\begin{proof}
Cf.\,proof of \cite[Theorem 4.2 and Proposition 5.1]{girsanovmarkov}.
\end{proof}
We can now apply \ref{girsanov} to build our favourite Markov process with symmetric rates out of $\alpha$:
\begin{cor}\label{corgir}
Let $\tilde{\alpha}$ be a Markov process on $\{0,1\}$ with generator
\begin{align}\label{deftildeq}
\tilde{Q}f(i):=\sqrt{s_0s_1}(f(1-i)-f(i))
\end{align}
for $f\colon\{0,1\}\to\mathbb{R}$, and write $\mathbb{P}^{\tilde{\alpha}}$ for its distribution with start in state $1$. Further, denote by $\mathbb{P}^{\alpha}$ the distribution of the Markov chain $\alpha$ with generator defined in \eqref{Q} and start in $1$. Then,
\begin{align}\label{radon}
\frac{\d\mathbb{P}^{\alpha}}{\d\mathbb{P}^{\tilde{\alpha}}}|_{\mathcal{F}_t}=\frac{h(\tilde{\alpha}(t))}{\sqrt{s_0}}\exp\left(\sqrt{s_0s_1}t-s_0\tilde{L}_t(0)-s_1\tilde{L}_t(1)\right)
\end{align}
where we wrote $\tilde{L}_t(i)=\int_0^t\delta_i(\tilde{\alpha}(s))\,\d s$ for the local times of $\tilde{\alpha}$ in state $i\in\{0,1\}$ up to time $t$.
\end{cor}
\begin{proof}[Proof]
Define $h:\{0,1\}\to\mathbb{R}$ as $h(0)=\sqrt{s_1}$ and $h(1)=\sqrt{s_0}$. Then,
\begin{align*}
\frac{\d\mathbb{P}^{\alpha}}{\d\mathbb{P}^{\tilde{\alpha}}}|_{\mathcal{F}_t}=&\frac{h(\tilde{\alpha}(t))}{h(\tilde{\alpha}(0))}\exp\left(-\int_0^t\frac{\tilde{Q}h(\tilde{\alpha}(s))}{h(\tilde{\alpha}(s))}\,\d s\right)
\\=&\frac{h(\tilde{\alpha}(t))}{\sqrt{s_0}}\exp\left(-\tilde{L}_t(0)\frac{\sqrt{s_0s_1}(\sqrt{s_0}-\sqrt{s_1})}{\sqrt{s_1}}-\tilde{L}_t(1)\frac{\sqrt{s_0s_1}(\sqrt{s_1}-\sqrt{s_0})}{\sqrt{s_0}}\right)
\\=&\frac{h(\tilde{\alpha}(t))}{\sqrt{s_0}}\exp\left(-\tilde{L}_t(0)(s_0-\sqrt{s_0s_1})-\tilde{L}_t(1)(s_1-\sqrt{s_0s_1})\right)
\\=&\frac{h(\tilde{\alpha}(t))}{\sqrt{s_0}}\exp\left(\sqrt{s_0s_1}t-s_0\tilde{L}_t(0)-s_1\tilde{L}_t(1)\right).
\end{align*}
\end{proof}
In particular, if $(Z,\alpha)$ is the Markov process with generator defined in \eqref{x-y} and $(\tilde{Z},\tilde{\alpha})$ another Markov process with symmetric generator
\begin{align*}
\tilde{L}f(x,i):=(i\kappa+\rho)\sum_{y\sim x}(f(y,i)-f(x,i))+\sqrt{s_0s_1}(f(x,1-i)-f(x,i))
\end{align*}
for test functions $f\colon\mathbb{Z}^d\times\{0,1\}\to\mathbb{R}$ and if we write $\mathbb{P}^{(Z,\alpha)}_{(0,1)}$ resp.\,$\tilde{\mathbb{P}}^{(Z,\tilde{\alpha})}_{(0,1)}$ for the distribution of $(Z,\alpha)$ resp.\,$(\tilde{Z},\tilde{\alpha})$ with start in $(0,1)$, then
\begin{align*}
\frac{\d\mathbb{P}^{(Z,\alpha)}_{(0,1)}}{\d\tilde{\mathbb{P}}^{(Z,\tilde{\alpha})}_{(0,1)}}|_{\mathcal{F}_t}=\frac{h(\tilde{\alpha}(t))}{\sqrt{s_0}}\exp\left(\sqrt{s_0s_1}t-s_0\tilde{L}_t(0)-s_1\tilde{L}_t(1)\right),
\end{align*}
since the generator of $Z$, conditioned on $\alpha$, equals that of $\tilde{Z}$, conditioned on $\tilde{\alpha}$, and since $\alpha$ resp.\,$\tilde{\alpha}$ is independent of $Z$ resp.\,$\tilde{Z}$.
\section{Proof of Theorem \ref{annasy1}}
This section is dedicated to the proof of Theorem \ref{annasy1}. For this, let $\xi=(\xi(x))_{x\in\mathbb{Z}^d}$ be a static field built out of Bernoulli distributed particles, i.\,e.\,, for each $x\in\mathbb{Z}^d$ the random variable $\xi(x)$ is independent and Bernoulli distributed with $\mathbb{P}(\xi(x)=1)=p=1-\mathbb{P}(\xi(x)=0)$. As the random environment is static, we only have to average over the movement of the switching random walk $X$ and the initial distribution of the Bernoulli field in order to determine the long-term behaviour of $\left<U(t)\right>$. Thus, the proof of Theorem \ref{annasy1} is based on a \emph{time-change} combined with existing results regarding the behaviour of random walks (without switching) in a Bernoulli field of particles. More precisely, let $\tilde{X}$ be a simple symmetric random walk without switching and with generator $\kappa\Delta$. Then it is well-known that, conditioned on $(\alpha(s))_{s\leq t}$, the endpoint $X(t)$ is equal to $\tilde{X}(L_t(1))$ in distribution. This will be the starting point of the following proof:
\begin{proof}[Proof of Theorem \ref{annasy1}] For $(x,i)\in\mathbb{Z}^d\times\{0,1\}$ let
\begin{align*}
\ell_t(x,i):=\int_0^t\delta_{(x,i)}(X(s),\alpha(s))\d s
\end{align*}
denote the local time of the process $(X,\alpha)$ in $(x,i)$ up to time $t\geq 0$. Then, for an arbitrary $\gamma\in[-\infty, \infty)$, we can rewrite $\left<U(t)\right>$, using the independence of the Bernoulli distribution in each spatial point $x\in\mathbb{Z}^d$, as
\begin{align*}
\left<U(t)\right>=&\left<\mathbb{E}_{(0,1)}^{(X,\alpha)}\left[\exp\left(\sum_{x\in\mathbb{Z}^d}\sum_{i\in\{0,1\}}\gamma\cdot i\cdot \xi(x)\ell_t(x,i)\right)\right]\right>
\\=&\mathbb{E}_{(0,1)}^{(X,\alpha)}\left[\prod_{x\in\mathbb{Z}^d}\left<\exp\left(\sum_{i\in\{0,1\}}\gamma\cdot i\cdot \xi(x)\ell_t(x,i)\right)\right>\right]
\\=&\mathbb{E}_{(0,1)}^{(X,\alpha)}\left[\prod_{x\in\mathbb{Z}^d}\left(p\e^{\gamma\ell_t(x,1)}+1-p\right)\right]
\\=&\mathbb{E}_{(0,1)}^{(X,\alpha)}\left[\exp\left(\sum_{x\in\mathbb{Z}^d}\log\left(p\e^{\gamma\ell_t(x,1)}+1-p\right)\mathds{1}_{\{\ell_t(x,1)>0\}}\right)\right].
\end{align*}
Now, let $\gamma=-\infty$. Then the annealed survival probability up to time $t$ reads
\begin{align*}
\left<U(t)\right>=&\mathbb{E}_{(0,1)}^{(X,\alpha)}\left[\exp\left(\sum_{x\in\mathbb{Z}^d}\log\left(1-p\right)\mathds{1}_{\{\ell_t(x,1)>0\}}\right)\right]
\\&=\mathbb{E}_{(0,1)}^{(X,\alpha)}\left[(1-p)^{\sum_{x\in\mathbb{Z}^d}\mathds{1}_{\{\ell_t(x,1)>0\}}}\right].
\end{align*}
Recall, that the walk $X$ can only move in the active state $1$ such that each newly visited point $x\in\mathbb{Z}^d$ is crossed by $X$ for the first time in state $1$. Therefore, 
\begin{align*}
\sum_{x\in\mathbb{Z}^d}\mathds{1}_{\{\ell_t(x,1)>0\}}=\sum_{x\in\mathbb{Z}^d}\mathds{1}_{\{\bar{\ell}_t(x)>0\}}=\colon R(t),
\end{align*}
where we write $\bar{\ell}_t(x)$ for the projection of $\ell_t(x,i)$, $i\in\{0,1\}$, on the first component and denote by $R(t)$ the range of the random walk $X$ up to time $t$, i.e.\,, the number of all distinct visited points up to time $t$ by $X$.  Changing the order of the expectations yields
\begin{align}\label{startlower}
\left<U(t)\right>=\mathbb{E}_{1}^{\alpha}\mathbb{E}_0^X\left[(1-p)^{R(t)}\right]=\mathbb{E}_{1}^{\alpha}\mathbb{E}_0^{\tilde{X}}\left[(1-p)^{\tilde{R}(L_t(1))}\right],
\end{align}
where $\mathbb{E}_0^{\tilde{X}}$ denotes the expectation with respect to $\tilde{X}$ and $\tilde{R}(t)$ the range of $\tilde{X}$ up to time $t$. For a fixed realization of $\alpha$, the inner expectation is nothing but the survival probability of the simple random walk $\tilde{X}$ among Bernoulli traps, to which we would like to apply the well-known result from \cite{antal} asserting that
\begin{align}\label{formulaantal}
\mathbb{E}_0^{\tilde{X}}\left[(1-p)^{\tilde{R}(t)}\right]=\exp\left(-c_d\kappa^{\frac{d}{d+2}}|\log(1-p)|^{\frac{2}{d+2}}t^{\frac{d}{d+2}}(1+o(1))\right)
\end{align}
as $t\to\infty$, where the constant $c_d$ is given by
\begin{align*}
c_d:=(d+2)d^{\frac{2}{d+2}-1}\lambda_d^{\frac{d}{d+2}}
\end{align*}
and $\lambda_d$ denotes the principal Dirichlet eigenvalue of $-\Delta$ on $[-1,1]^d\subseteq\mathbb{R}^d$. Heuristically, we would expect to have
\begin{align*}
\mathbb{E}_0^{\tilde{X}}\left[(1-p)^{\tilde{R}(L_t(1))}\right]=\exp\left(-c_d\kappa^{\frac{d}{d+2}}(\log(1-p))^{\frac{2}{d+2}}L_t(1)^{\frac{d}{d+2}}(1+o(1))\right)
\end{align*}
as $t\to\infty$ and then apply Lemma \ref{varadhanversion} to finish the proof of part (a). However, we have to be careful with the error term, as $L_t(1)$ is not a deterministic time. To handle this, note that for any $\delta>0$ and using the Markov property,
\begin{align*}
\left<U(t)\right>&\geq \mathbb{E}_1^\alpha\left[\int_0^1\mathbb{E}_0^{\tilde{X}}\left[(1-p)^{\tilde{R}((a+\delta)t)}\right]1_{\{L_t(1)\in t B_{\delta}(a)\}}\d a\right]
\\&\geq\mathbb{E}_1^\alpha\left[\int_0^1\mathbb{E}_0^{\tilde{X}}\left[(1-p)^{\tilde{R}((a-\delta)t)}\right]\mathbb{E}_0^{\tilde{X}}\left[(1-p)^{\tilde{R}(2\delta t)}\right]1_{\{L_t(1)\in t B_{\delta}(a)\}}\d a\right]
\\&=\mathbb{E}_1^{\alpha}\left[\int_0^1\exp\left(-K\left(((a-\delta)t)^{\frac{d}{d+2}}+(2\delta t)^{\frac{d}{d+2}}\right)(1+o(1))\right)1_{\{L_t(1)\in t B_{\delta}(a)\}}\d a\right],
\end{align*}
as $t\to\infty$ and with $K:=c_d\kappa^{\frac{d}{d+2}}|\log(1-p)|^{\frac{2}{d+2}}$, where in the last step we used \eqref{formulaantal} and where the error in $o(1)$ does not depend on $\alpha$. Hence, as we condition on the event ${\{L_t(1)\in t B_{\delta}(a)\}}$ in the expectation and so $(a-\delta)t\geq L_t(1)-2\delta t$ in this event,
\begin{align*}
\left<U(t)\right>&\geq \mathbb{E}_1^{\alpha}\left[\int_0^1\exp\left(-K\left(L_t(1)^{\frac{d}{d+2}}+(2\delta t)^{\frac{d}{d+2}}\right)(1+o(1))\right)1_{\{L_t(1)\in t B_{\delta}(a)\}}\d a\right]
\\&\geq \exp\left(-K(2\delta t)^{\frac{d}{d+2}}\right)\mathbb{E}_1^{\alpha}\left[\exp\left(-K L_t(1)^{\frac{d}{d+2}}(1+o(1))\right)\right].
\end{align*}
Applying Lemma \eqref{varadhanversion} we obtain
\begin{align*}
\liminf_{t\to\infty}\frac{1}{t^{\frac{d}{d+2}}}\log\left<U(t)\right>\geq -K(2\delta)^{\frac{d}{d+2}}-K\left(\frac{s_0}{s_0+s_1}\right)^{\frac{d}{d+2}}
\end{align*}
for all $\delta>0$ such that we can let $\delta$ tend to $0$. For the upper bound, we proceed in a similar way and estimate
\begin{align*}
\left<U(t)\right>&\leq \mathbb{E}_1^\alpha\left[\int_0^1\mathbb{E}_0^{\tilde{X}}\left[(1-p)^{\tilde{R}((a-\delta)t)}\right]1_{\{L_t(1)\in t B_{\delta}(a)\}}\d a\right]
\\&\leq \mathbb{E}_1^{\alpha}\left[\int_0^1\exp\left(-K\left(((a-\delta)t)^{\frac{d}{d+2}}-(2\delta t)^{\frac{d}{d+2}}\right)(1+o(1))\right)1_{\{L_t(1)\in t B_{\delta}(a)\}}\d a\right]
\\&\leq \exp\left(K(2\delta t)^{\frac{d}{d+2}}\right)\mathbb{E}_1^{\alpha}\left[\exp\left(-K L_t(1)^{\frac{d}{d+2}}(1+o(1))\right)\right],
\end{align*}
since $(a-b)t\geq L_t(1)-2\delta t$ on the event $\{L_t(1)\in t B_{\delta}(a)\}$. Applying Lemma \eqref{varadhanversion} we again  obtain
\begin{align*}
\limsup_{t\to\infty}\frac{1}{t^{\frac{d}{d+2}}}\log\left<U(t)\right>\leq K(2\delta)^{\frac{d}{d+2}}-K\left(\frac{s_0}{s_0+s_1}\right)^{\frac{d}{d+2}}
\end{align*}
for all $\delta>0$. This finishes the proof of case part (a). 

Now, let $\gamma>0$. Then
\begin{align*}
\left<U(t)\right>&=\mathbb{E}_{(0,1)}^{(X,\alpha)}\left[\exp\left(\sum_{x\in\mathbb{Z}^d}\log\left(p\e^{\gamma\ell_t(x,1)}+1-p\right)\mathds{1}_{\{\ell_t(x,1)>0\}}\right)\right]
\\&=\mathbb{E}_{(0,1)}^{(X,\alpha)}\left[\exp\left(\sum_{x\in\mathbb{Z}^d}\left(\log(p)+\log\left(\e^{\gamma\ell_t(x,1)}+\frac{1-p}{p}\right)\right)\mathds{1}_{\{\ell_t(x,1)>0\}}\right)\right]
\\&=\mathbb{E}_{(0,1)}^{(X,\alpha)}\left[p^{R(t)}\cdot\exp\left(\sum_{x\in\mathbb{Z}^d}\log\left(\e^{\gamma\ell_t(x,1)}+\frac{1-p}{p}\right)\mathds{1}_{\{\ell_t(x,1)>0\}}\right)\right].
\end{align*}
Using the asymptotics \eqref{formulaantal} of the first term in the expectation (with $1-p$ replaced by $p$) and time-change again, we can lower-bound
\begin{align*}
\left<U(t)\right>\geq& \mathbb{E}_{(0,1)}^{(X,\alpha)}\left[p^{R(t)}\cdot\exp\left(\sum_{x\in\mathbb{Z}^d}\log\left(\e^{\gamma\ell_t(x,1)}\right)\mathds{1}_{\{\ell_t(x,1)>0\}}\right)\right]
\\=&\mathbb{E}_1^{\alpha}\mathbb{E}_0^{\tilde{X}}\left[p^{\tilde{R}(L_t(1))}\exp\left(\gamma L_t(1)\right)\right]
\\=&\mathbb{E}_1^{\alpha}\left[\exp\left(\gamma L_t(1)-c_d\kappa^{\frac{d}{d+2}}\log(p)^{\frac{2}{d+2}}L_t(1)^{\frac{d}{d+2}}(1+o(1))\right)\right], \quad t\to\infty,
\end{align*}
where the error term does not depend on the fluctuations of $\alpha$, which is seen in an analogous way as in the case $\gamma<0$. We omit repeating the details. Now, recall that by Corollary \ref{LDP}, the normalized local times $\small\frac{1}{t}L_t(1)$ of $\alpha$ in $1$ satisfy a large deviation principle with rate function $I$ defined in \eqref{I}. Using Varadhan's Lemma (cf.\,\cite[Theorem 3.3.1]{koenigldp}), we can deduce
\begin{align*}
\lim_{t\to\infty}\frac{1}{t}\log\left<U(t)\right>\geq \sup_{a\in[0,1]}\left\{\gamma a-I(a)\right\}=\gamma-\inf_{a\in[0,1]}\left\{(1-a)\gamma+I(a)\right\},
\end{align*}
as $t\to\infty$, since the term
\begin{align*}
L_t(1)^{\frac{d}{d+2}}=t^{-\frac{2}{d+2}}\left(\frac{1}{t}L_t(1)\right)^{\frac{d}{d+2}}
\end{align*}
vanishes on scale $t$ for $t\to\infty$. On the other hand, we can upper-bound
\begin{align*}
\left<U(t)\right>&=\mathbb{E}_{(0,1)}^{(X,\alpha)}\left[\exp\left(\sum_{x\in\mathbb{Z}^d}\log\left(p\e^{\gamma\ell_t(x,1)}+1-p\right)\mathds{1}_{\{\ell_t(x,1)>0\}}\right)\right]
\\&\leq\mathbb{E}_{(0,1)}^{(X,\alpha)}\left[\exp\left(\sum_{x\in\mathbb{Z}^d}\gamma\ell_t(x,1)\mathds{1}_{\{\ell_t(x,1)>0\}}\right)\right]=\mathbb{E}_{1}^{\alpha}\left[\exp(\gamma L_t(1))\right]
\end{align*}
and thus
\begin{align*}
\lim_{t\to\infty}\frac{1}{t}\log\left<U(t)\right>\leq \gamma-\inf_{a\in[0,1]}\left\{(1-a)\gamma+I(a)\right\}.
\end{align*}
For an explicit expression, we calculate the minimizer of the function $f(a)=(1-a)\gamma+I(a)$ and find that
\begin{align*}
\gamma-\inf_{a\in[0,1]}\left\{(1-a)\gamma+I(a)\right\}=\gamma-s_1-\frac{(\gamma+s_0-s_1)^2-s_0s_1}{\sqrt{\gamma^2+2\gamma(s_0-s_1)+(s_0+s_1)^2}}.
\end{align*}
\end{proof}
\section{Proof of Theorem \ref{annasy2}}
In this section, we prove Theorem \ref{annasy2}, in which the underlying environment $\xi$ is dynamic and consists of one single particle moving independently of $X$. More precisely, $\xi$ is the Markov process on $\{0,1\}^{\mathbb{Z}^d}$ given by
\begin{align*}
\xi(x,t)=\delta_x(Y(t)),
\end{align*}
where $Y=(Y(t))_{t\geq 0}$ is a continuous-time simple symmetric random walk on $\mathbb{Z}^d$ with total jump rate $2d\rho$ for a constant $\rho>0$ and starting in the origin. Hence,
\begin{align}\label{defuorigin}
\left<U(t)\right>=\mathbb{E}_0^Y\mathbb{E}_{(0,1)}^{(X,\alpha)}\left[\exp\left(\gamma\int_0^t\delta_{(0,1)}(X(s)-Y(s),\alpha(s))\,\dt d s\right)\right],
\end{align}
where $\mathbb{E}_0^Y$ denotes the expectation with respect to $Y$. Set $Z:=X-Y$. Then $(Z,\alpha)$ has the generator
\begin{align*}
\bar{\mathcal{L}}f(z,i)=\sum_{y\sim z}(i\kappa+\rho)(f(y,i)-f(z,i))+s_i(f(z,1-i)-f(z,i))
\end{align*} 
for $z\in\mathbb{Z}^d$, $i,j\in\{0,1\}$ and a test function $f:\mathbb{Z}^d\times \{0,1\}\to\mathbb{R}$, and thus
\begin{align}\label{onesurv}
\left<U(t)\right>=\mathbb{E}_{(0,1)}^{(Z,\alpha)}\left[\exp\left(\gamma\int_0^t\delta_{(0,1)}(Z(s),\alpha(s))\,\dt d s\right)\right]=:v(0,1,t).
\end{align}
Especially, the new potential $\tilde\xi(z,i):=\delta_{(0,1)}(x,i)$ does not depend on the time any more. Using the Feynman-Kac formula, we further see that \eqref{onesurv} is the solution to
\begin{align}\label{defv}
\left\{\begin{array}{lllr}
\frac{\dt d}{\dt dt} v(y,i,t)&=&(i\kappa+\rho)\Delta v(y,i, t)+Q v(y,i,t)+\gamma\cdot \delta_{(0,1)}(y,i)v(y,i,t), &t>0
\\[8pt]v(y,i,0)&=&i,
\end{array}\right.
\end{align}
with $(y,i)=(0,1)$. In the following, we shall use either of the representations \eqref{defuorigin}, \eqref{onesurv} or \eqref{defv}, depending on what is to be proven. We start with the proof of theorem \ref{annasy2}(a) and show theorem \ref{annasy2}(b) separately, as different methods are used in the case $\gamma<0$ and $\gamma>0$ respectively.
\begin{proof}[Proof of Theorem \ref{annasy2}(a)] Let $\gamma\in(-\infty,0)$ and denote by $p(y,i,t)$ the probability mass function of the switching diffusion $(Z,\alpha)$ with start in $(0,1)$. Then we get the representation
\begin{align}\label{u}
v(0,1,t)=1+\gamma\int_0^t p(0,1,s)v(0,1,t-s)\,\d s,
\end{align}
to which we want to apply the Laplace transform. Denoting by $\hat{v_1}(\lambda)$ resp. $\hat{p_1}(\lambda)$ the Laplace transform of $v(0,1,t)$ resp. $p(0,1,t)$, and solving \eqref{u} for $\hat{v_1}(\lambda)$, we arrive at
\begin{align}\label{uhat}
\hat{v_1}(\lambda)=\frac{1}{\lambda}\cdot\frac{1}{1-\gamma\hat{p_1}(\lambda)}.
\end{align}
Our next aim is therefore to compute $\hat{p_1}(\lambda)$. For this, note that the probability density function $p$ satisfies the system of equations 
\begin{align}\label{solut}
\left\{\begin{array}{lllr}\frac{\d}{\d t}p(y,1,t)&=&(\kappa+\rho)\Delta p(y,1,t)+s_0p(y,0,t)-s_1p(y,1,t)), &t>0,
\\[8pt]\frac{\d}{\d t}p(y,0,t)&=&\rho\Delta p(y,0,t)+s_1p(y,1,t)-s_0p(y,0,t), &t>0,
\\[8pt]p(y,i,0)&=&\delta_{(0,1)}(y,i).
\end{array}\right.
\end{align}
Forth order systems of the form \eqref{newsystem} with two different diffusion constants have been studied in \cite{doubleI}. For the convenience of the reader, we will include the first steps and ideas to calculate the solution of \eqref{solut}. We denote by $\hat{p}_i(y,\lambda)$ the Laplace transform of $p(y,i,\cdot)$ and apply this to \eqref{solut}, using the initial condition, to obtain the new system
\begin{align}
0=&(\kappa+\rho)\Delta\hat{p}_1(y,\lambda)-(\lambda+s_1)\hat{p}_1(y,\lambda)+s_0\hat{p}_0(y,\lambda)+\delta_0(y),\label{1}
\\[8pt]0=&\rho\Delta\hat{p}_0(y,\lambda)-(\lambda+s_0)\hat{p}_0(y,\lambda)+s_1\hat{p}_1(y,\lambda)\label{2},
\end{align}
which, after solving \eqref{2} for $\hat{p}_1(y,\lambda)$ and applying this to \eqref{1}, translates in to the forth-order equation
\begin{align}\label{newsystem}
\left(\Delta^2-\left(\frac{s_1+\lambda}{\kappa+\rho}+\frac{s_0+\lambda}{\rho}\right)\Delta+\frac{(s_1+\lambda)(s_0+\lambda)-s_0s_1}{(\kappa+\rho)\rho}\right)\hat{p}_0(y,\lambda)=\frac{s_1}{(\kappa+\rho)\rho}\delta_0(y).
\end{align}
Forth order systems of the form \eqref{newsystem} are known to have the solution
\begin{align*}
\hat{p}_0(y,\lambda)=\frac{s_1}{2(\kappa+\rho)\rho (a^2-b^2)}\left(\frac{1}{a}\e^{a|y|}-\frac{1}{b}\e^{b|y|}\right),
\end{align*}
in dimension $d=1$, where
\begin{align*}
a,b=-\frac{\sqrt{\lambda(\kappa+\rho)+s_1\rho+s_0(\kappa+\rho)\pm\sqrt{\kappa^2+2\lambda\kappa(s_1\rho-s_0(\kappa+\rho))+(s_1\rho+s_0(\kappa+\rho))^2}}}{\sqrt{2\rho(\kappa+\rho)}}.
\end{align*}
Using the relation between $\hat{p}_1(y,\lambda)$ and $\hat{p}_0(y,\lambda)$ and inserting $y=0$ we obtain
\begin{align*}
\hat{p}_1(0,\lambda)= \frac{-(\lambda+s_0+\rho ab)}{2\rho(\kappa+\rho)ab(a+b)}\sim \frac{1}{\sqrt{\lambda}}\cdot\frac{s_0}{2\sqrt{(s_0+s_1)(s_0(\rho+\kappa)+s_1\rho)}}, \quad \lambda\to 0,
\end{align*}
as long and tedious calculations show. In dimension $d=2$ we proceed in a similar way to find
\begin{align*}
\hat{p}_1(0,\lambda)\sim \frac{s_0}{4\pi(s_1\rho+s_0(\kappa+\rho))}\log\left(\frac{1}{\lambda}\right),\quad \lambda\to 0.
\end{align*}
Thus, we deduce from \eqref{uhat} that
\begin{align*}
\hat{v}_1(\lambda)\sim \left\{\begin{array}{ll}\displaystyle\frac{1}{\sqrt{\lambda}}\frac{2\sqrt{(s_0+s_1)(s_0(\rho+\kappa)+s_1\rho)}}{s_0|\gamma|}, &d=1,
\\[12pt]\displaystyle\frac{4\pi(s_1\rho+s_0(\kappa+\rho))}{s_0|\gamma|\lambda\log\left(\frac{1}{\lambda}\right)}, &d=2,\end{array}\right.
\end{align*} 
as $\lambda\to 0$. 
Using the Littlewood-Hardy Tauberian theorem we finally arrive at
\begin{align*}
v(1,0,t)\sim\left\{\begin{array}{ll} \displaystyle\frac{2\sqrt{(s_0+s_1)(s_0(\rho+\kappa)+s_1\rho)}}{\sqrt{\pi}s_0|\gamma|}\frac{1}{\sqrt{t}}, &d=1,\\[12pt]\displaystyle\frac{4\pi(s_1\rho+s_0(\kappa+\rho))}{s_0|\gamma|\log(t)}, &d=2,\end{array}\right.
\end{align*} 
as $\lambda\to\infty$. Next, let $d\geq 3$ and denote by
\begin{align*}
G_d(x,i):=\int_0^{\infty}p_d(x,i,t)\,\d t
\end{align*}
the Green's function of $(Z,\alpha)$ in $(x,i)$, which has the probabilistic representation
\begin{align*}
G_d(x,i)=\mathbb{E}_{(x,i)}^{(Z,\alpha)}\left[\int_0^{\infty}\delta_{(0,1)}(Z(s),\alpha(s))\,\d s\right].
\end{align*}
Hence, for all $(x,i)\in\mathbb{Z}^d\times\{0,1\}$ the quantity
\begin{align*}
v(x,i):=\lim_{t\to\infty}v(x,i,t)=\mathbb{E}_{(x,i)}^{(Z,\alpha)}\left[\exp\left(\gamma\int_0^{\infty}\delta_{(0,1)}(Z(s),\alpha(s))\,\d s\right)\right]
\end{align*}
lies in $(0,1)$. Moreover, $v$ solves the boundary value problem
\begin{align*}
\left\{\begin{array}{llll}0&=&(i\kappa+\rho)\Delta v(x,i)+\gamma\delta_{(0,1)}(x,i)v(x,i), &(x,i)\in\mathbb{Z}^d\times\{0,1\},\\[8pt]1&=&\lim_{\|x\|\to\infty}v(x,i), &i\in\{0,1\},\end{array}\right.
\end{align*}
and can therefore be written as
\begin{align*}
v(0,1)=1+\gamma\int_0^{\infty}p_d(0,1,t)v(0,1)\,\d t=1+\gamma v(0,1)G_d(0,1)
\end{align*}
in point $(0,1)$. Solving for $v(0,1)$ gives
\begin{align*}
v(0,1)=\frac{1}{1-\gamma G_d(0,1)}.
\end{align*}
The survival probability converges therefore to a non-trivial limit in $(0,1)$ in all dimensions $d\geq 3$.
\end{proof}
We now continue with the case $\gamma>0$ of catalysts. Recall the two-state Markov chain $\tilde{\alpha}$ with symmetric generator \eqref{deftildeq}. Before proving theorem \ref{annasy2}(b), we need two statements that are highly inspired by \cite[Lemma 2.2 and Lemma 2.3]{movingcat}:
\begin{lemma}\label{lemmaab}
Let $r(t)=t\log^2(t)$, $Q_{r(t)}=[-r(t),r(t)]^d\cap\mathbb{Z}^d$ and $V:\mathbb{Z}^d\times\{0,1\}\to\mathbb{R}$ a bounded function. Further, abbreviate
\begin{align*}
A_t:=\int_0^tV(Z(s),\tilde{\alpha}(s))\,\d s.
\end{align*}
Then, the following holds true:
\begin{enumerate}
\item[(a)] As $t\to\infty$,
\begin{align*}
\mathbb{E}^{(Z,\tilde{\alpha})}_{(0,1)}&\left[\e^{A_t}\right]=(1+o(1))\sum_{z\in\mathbb{Z}^d}\mathbb{E}_{(0,1)}^{(X,\tilde{\alpha})}\mathbb{E}_z^Y\left[\e^{A_t}\cdot\delta_0(Y(t))\mathds{1}_{\{X(t)\in Q_{r(t)}\}}\right].
\end{align*}
\item[(b)] As $t\to\infty$,
\begin{align*}
\sum_{y\in Q_{r(t)}}\mathbb{E}_{(0,1)}^{(X,\tilde{\alpha})}&\mathbb{E}_0^Y\left[\e^{A_t}\cdot\delta_y(X(t))\delta_y(Y(t))\right]=(1+o(1))\sum_{y\in\mathbb{Z}^d}\mathbb{E}_{(0,1)}^{(X,\tilde{\alpha})}\mathbb{E}_0^Y\left[\e^{A_t}\cdot\delta_y(X(t))\delta_y(Y(t))\right].
\end{align*}
\end{enumerate}
\end{lemma}
\begin{proof}
Note that $A_t\in[0, Mt]$ with $M:=\sup_{(x,i)\in\mathbb{Z}^d\times\{0,1\}}V(x,i)$. Then,
\begin{align*}
\mathbb{E}^{(Z,\tilde{\alpha})}_{(0,1)}\left[\e^{A_t}\right]=\sum_{z\in\mathbb{Z}^d}\mathbb{E}_{(0,1)}^{(X,\tilde{\alpha})}\mathbb{E}_0^Y\left[\e^{A_t}\delta_z(Y(t))\right]=\sum_{z\in\mathbb{Z}^d}\mathbb{E}_{(0,1)}^{(X,\tilde{\alpha})}\mathbb{E}_z^Y\left[\e^{A_t}\delta_0(Y(t))\right]
\end{align*}
using Fubini and a time reversal for $Y$. 
In order to prove part (a) we have to check that
\begin{align*}
a(t):=\frac{\sum_{z\in\mathbb{Z}^d}\mathbb{E}_{(0,1)}^{(X,\tilde{\alpha})}\mathbb{E}_z^Y\left[\e^{A_t}\delta_0(Y(t))\right]-\sum_{z\in Q_{r(t)}}\mathbb{E}_{(0,1)}^{(X,\tilde{\alpha})}\mathbb{E}_z^Y\left[\e^{A_t}\delta_0(Y(t))\mathds{1}_{\{X(t)\in Q_{r(t)}\}}\right]}{\sum_{z\in\mathbb{Z}^d}\mathbb{E}_{(0,1)}^{(X,\tilde{\alpha})}\mathbb{E}_z^Y\left[\e^{A_t}\delta_0(Y(t))\right]},
\end{align*}
converges to zero as $t\to\infty$, which is done in a similar way as in the proofs of \cite[Lemma 2.2]{movingcat}, such that we only highlight the differences. Splitting $\mathbb{Z}^d$ in $Q_{r(t)}$ and it's complement, upper-bounding $\e^{A_t}$ by $\e^{M t}$ and using a time reversal for $Y$ again, we obtain the bound
\begin{align}\label{a(t)}
a(t)\leq& \frac{\e^{M t}\left(\mathbb{P}_0^Y(Y(t)\notin Q_{r(t)})+\mathbb{P}_{(0,1)}^{(X,\tilde{\alpha})}(X(t)\notin Q_{r(t)})\right)}{\mathbb{P}_0^Y(Y(t)=0)}\nonumber
\\=&\frac{\e^{M t}\left(\mathbb{P}_0^Y(Y(t)\notin Q_{r(t)})+\mathbb{P}_1^{\tilde{\alpha}}\mathbb{P}_{0}^{\tilde{X}}(\tilde{X}(\tilde{L}_t(1))\notin Q_{r(t)}\right)}{\mathbb{P}_0^Y(Y(t)=0)},
\end{align}
where $\tilde{X}$ is a simple symmetric random walk without switching and with generator $\kappa\Delta$ and $\tilde{L}_t(1)$ denotes the local time of $\tilde{\alpha}$ in state $1$ up to time $t$. Now, \cite[Lemma 4.3]{garmol} asserts that
\begin{align*}
\mathbb{P}_0^Y(Y(t)\notin Q_{r(t)})\leq 2^{d+1}\e^{-r(t)\log\left(\frac{r(t)}{d\rho t}\right)+r(t)}
\end{align*} 
such that for our choice of $r(t)$ and for sufficiently large $t$,
\begin{align*}
\mathbb{P}_0^Y(Y(t)\notin Q_{r(t)})\leq \e^{-r(t)},
\end{align*}
as a quick estimation shows. Analogously,
\begin{align*}
\mathbb{P}_1^{\tilde{\alpha}}\mathbb{P}_{0}^{\tilde{X}}\left(\tilde{X}(L_t(1))\notin Q_{r(t)}\right)\leq& \mathbb{P}_1^{\tilde{\alpha}}\left[2^{d+1}\exp\left(-r(t)\log\left(\frac{r(t)}{d\kappa \tilde{L}_t(1)}\right)+r(t)\right)\right]
\\=&2^{d+1}\exp\left(-r(t)\left(\log\left(\frac{r(t)}{d\kappa}\right)-1\right)\right)\mathbb{E}_1^{\tilde{\alpha}}\left[ \exp\left(r(t)\log(\tilde{L}_t(1))\right)\right]
\\\leq&2^{d+1}\exp\left(-r(t)\log\left(\frac{r(t)}{d\kappa t}\right)+r(t)\right).
\end{align*}
Thus, we have again
\begin{align*}
\mathbb{P}_1^{\tilde{\alpha}}\mathbb{P}_{0}^{\tilde{X}}\left(\tilde{X}(\tilde{L}_t(1))\notin Q_{r(t)}\right)\leq\e^{-r(t)}
\end{align*}
for sufficiently large $t$. This shows that the numerator of \eqref{a(t)} converges exponentially in $t$ to zero, whereas its denominator converges only polynomially. Hence, $a(t)\to 0$ for $t\to\infty$. 

In order to prove part (b), we define
\begin{align*}
b(t):=\frac{\sum_{y\notin Q_{r(t)}}\mathbb{E}_{(0,1)}^{(X,\tilde{\alpha})}\mathbb{E}_0^Y\left[\e^{A_t}\delta_y(X(t))\delta_y(Y(t))\right]}{\sum_{y\in \mathbb{Z}^d}\mathbb{E}_{(0,1)}^{(X,\tilde{\alpha})}\mathbb{E}_0^Y\left[\e^{A_t}\delta_y(X(t))\delta_y(Y(t))\right]}
\end{align*}
and proceed in a similar way to obtain the upper bound
\begin{align}\label{b(t)}
b(t)\leq\frac{\e^{M t}\left(\mathbb{P}_{(0,1)}^{(X,\tilde{\alpha})}(X(t)\notin Q_{r(t)})\mathbb{P}_0^Y(Y(t)\notin Q_{r(t)})\right)}{\mathbb{P}_{(0,1)}^{(X,\tilde{\alpha})}(X(t)=0)\mathbb{P}_0^Y(Y(t)=0)}.
\end{align}
From the proof of part (a) we already know that the nominator decays exponentially in $t$ as $t\to\infty$. Moreover,
\begin{align*}
\mathbb{P}_{(0,1)}^{(X,\tilde{\alpha})}(X(t)=0)=\mathbb{P}_1^{\tilde{\alpha}}\mathbb{P}_0^{\tilde{X}}(\tilde{X}(\tilde{L}_t(1))=0),
\end{align*}
where $\mathbb{P}_0^{\tilde{X}}(\tilde{X}(\tilde{L}_t(1))=0)$ decays polynomially in $\tilde{L}_t(1)$. Hence, lemma \ref{varadhanversion} asserts that
\begin{align*}
\mathbb{P}_1^{\tilde{\alpha}}\mathbb{P}_0^{\tilde{X}}(\tilde{X}(\tilde{L}_t(1))=0)=\mathbb{P}_0^{\tilde{X}}\left(\tilde{X}\left(\frac{s_0}{s_0+s_1}t\right)=0\right)
\end{align*}
for sufficiently large $t$. Thus, the denominator of the right hand-side of \eqref{b(t)} decays polynomially in $t$ and therefore $b(t)\to 0$ as $t\to\infty$. 
\end{proof}
We are now ready to prove theorem \ref{annasy2}(b).
\begin{proof}[Proof of Theorem \ref{annasy2}(b)]
Let $\tilde{\alpha}$ be the two-state Markov chain with symmetric generator $\tilde{Q}$ defined in \eqref{deftildeq} and denote by $\ell_t(x,i)$ resp.\,$\tilde{\ell}_t(x,i)$ the local times of $(Z,\alpha)$ resp.\,$(Z,\tilde{\alpha})$ in $(x,i)$ up to time $t$. Then, combining the representation \eqref{onesurv} with Corollary \ref{corgir}, the annealed number of individuals up to time $t$ reads
\begin{align}\label{dar}
\left<U(t)\right>=&\mathbb{E}^{(Z,\alpha)}_{(0,1)}\left[\exp\left(\gamma\int_0^t\delta_{(0,1)}(Z(s),\alpha(s))\,\d s\right)\right]\nonumber
\\=&\mathbb{E}^{(Z,\alpha)}_{(0,1)}\left[\exp\left(\gamma\ell_t(0,1)\right)\right]\nonumber
\\=&\mathbb{E}^{(Z,\tilde{\alpha})}_{(0,1)}\left[\exp\left(\sqrt{s_0s_1}t-s_0\tilde{L}_t(0)-s_1\tilde{L}_t(1)+\gamma\tilde{\ell}_t(0,1)\right)\right]\nonumber
\\=&\mathbb{E}^{(Z,\tilde{\alpha})}_{(0,1)}\left[\exp\left(\sqrt{s_0s_1}t-s_0\sum_{y\in\mathbb{Z}^d}\tilde{\ell}_t(y,0)-s_1\sum_{y\in\mathbb{Z}^d}\tilde{\ell}_t(y,1)+\gamma\tilde{\ell}_t(0,1)\right)\right]\nonumber
\\=&\e^{\sqrt{s_0s_1}t}\mathbb{E}^{(Z,\tilde{\alpha})}_{(0,1)}\left[\exp\left(\int_0^tV(Z(s),\tilde{\alpha}(s))\,\d s\right)\right]
\end{align}
where we define
\begin{align*}
V(z,i):=-s_0\delta_0(i)-s_1\delta_{1}(i)+\gamma\delta_{(0,1)}(z,i)
\end{align*}
for $(z,i)\in\mathbb{Z}^d\times\{0,1\}$. 
Let us start with the  upper bound, which is done in a similar way as in the proof of \cite[Theorem 1.2]{movingcat}. Applying Lemma \ref{lemmaab}(a) yields
\begin{align}\label{upperb}
\mathbb{E}^{(Z,\tilde{\alpha})}_{(0,1)}&\left[\exp\left(\int_0^tV(Z(s),\tilde{\alpha}(s))\,\d s\right)\right]\nonumber
\\=&(1+o(1))\sum_{z\in\mathbb{Z}^d}\mathbb{E}_{(0,1)}^{(X,\tilde{\alpha})}\mathbb{E}_z^Y\left[\exp\left(\int_0^tV(X(s)-Y(s),\tilde{\alpha}(s))\,\d s\right)\delta_0(Y(t))\mathds{1}_{\{X(t)\in Q_{r(t)}\}}\right]\nonumber
\\\leq& (1+o(1))\sum_{z\in\mathbb{Z}^d}\mathbb{E}_{(0,1)}^{(X,\tilde{\alpha})}\mathbb{E}_z^Y\left[\exp\left(\int_0^tV(X(s)-Y(s),\tilde{\alpha}(s))\,\d s\right)\mathds{1}_{\{X(t)-Y(t)\in Q_{r(t)}\}}\right]\nonumber
\\=&(1+o(1))\sum_{z\in\mathbb{Z}^d}\mathbb{E}_{(z,1)}^{(Z,\tilde{\alpha})}\left[\exp\left(\int_0^tV(Z(s),\tilde{\alpha}(s))\,\d s\right)\mathds{1}_{\{Z(t)\in Q_{r(t)}\}}\right].
\end{align}
Denote with $\left(\cdot,\cdot\right)$ the inner product in $\ell^2(\mathbb{Z}^d\times\{0,1\})$ with corresponding norm $\|\cdot\|$ and let 
\begin{align*}
\lambda:=\sup\text{Sp}(\tilde{L}+V)
\end{align*}
be the largest eigenvalue of the bounded and self-adjoint operator $\tilde{L}+V$. Then, applying the spectral representation to the right hand-side of \eqref{upperb} and proceeding in the standard way we obtain the upper bound
\begin{align*}
\mathbb{E}^{(Z,\tilde{\alpha})}_{(0,1)}\left[\exp\left(\int_0^tV(Z(s),\tilde{\alpha}(s))\d s\right)\right]\leq& (1+o(1))\left(\e^{(\tilde{L}+V)t}\mathds{1}_{Q_{r(t)}},\mathds{1}_{Q_{r(t)}}\right)
\\\leq&(1+o(1))\e^{t\lambda}\|\mathds{1}_{Q_{r(t)}}\|^2
\\\leq&(1+o(1))\e^{t\lambda}|Q_{r(t)}|
\\=&(1+o(1))\e^{t\lambda}(2t\log^2(t))^d.
\end{align*}
As $(2t\log^2(t))^d$ grows only polynomially, we have
\begin{align*}
\lim_{t\to\infty}\frac{1}{t}\log\mathbb{E}^{(Z,\tilde{\alpha})}_{(0,1)}\left[\exp\left(\int_0^tV(Z(s),\tilde{\alpha}(s))\,\d s\right)\right]\leq \lambda.
\end{align*}
For the lower bound, we proceed as proof of \cite[Theorem 1.2]{movingcat} to obtain
\begin{align*}
\mathbb{E}^{(Z,\tilde{\alpha})}_{(0,1)}\left[\exp\left(\int_0^tV(Z(s),\tilde{\alpha}(s))\d s\right)\right]\geq \frac{1}{|Q_{r(t)}|}\left(\sum_{y\in Q_{r(t)}}\mathbb{E}^{(X,\tilde{\alpha})}_{(0,1)}\mathbb{E}_0^Y\left[\e^{A_{t/2}}\delta_y(X(t/2)\delta_y(Y(t/2))\right]\right)^2
\end{align*}
for $A_t:=\exp\left(\int_0^t V(X(s)-Y(s),\tilde{\alpha}(s)\,\d s\right)$. Then, applying Lemma \ref{lemmaab}(b) yields
\begin{align*}
\mathbb{E}^{(Z,\tilde{\alpha})}_{(0,1)}\left[\exp\left(\int_0^tV(Z(s),\tilde{\alpha}(s))\,\d s\right)\right]\geq&\frac{(1+o(1))}{|Q_{r(t)}|}\left(\sum_{y\in \mathbb{Z}^d}\mathbb{E}^{(X,\tilde{\alpha})}_{(0,1)}\mathbb{E}_0^Y\left[\e^{A_{t/2}}\delta_y(X(t/2)\delta_y(Y(t/2))\right]\right)^2
\\=&\frac{(1+o(1))}{|Q_{r(t)}|}\left(\mathbb{E}^{(X,\tilde{\alpha})}_{(0,1)}\mathbb{E}_0^Y\left[\e^{A_{t/2}}\delta_0(X(t/2)-Y(t/2))\right]\right)^2
\\\geq &\frac{(1+o(1))}{|Q_{r(t)}|}\left(\e^{(\tilde{L}+V)\frac{t}{2}}\delta_0,\delta_0\right)^2.
\end{align*}
Now, we restrict the operator $\tilde{L}+V$ to finite boxes $Q_n:=([-n,n]^d\cap\mathbb{Z}^d)\times\{0,1\}$ and apply the Perron-Frobenius theorem for non-negative irreducible matrices to derive the existence of a largest eigenvalue $\lambda_n$ of $\tilde{L}+V$ on $Q_n:=([-n,n]^d\cap\mathbb{Z}^d)\times\{0,1\}$, for which
\begin{align*}
\lim_{t\to\infty}\frac{1}{t}\log\mathbb{E}^{(Z,\tilde{\alpha})}_{(0,1)}\left[\exp\left(\int_0^tV(Z(s),\tilde{\alpha}(s))\,\d s\right)\right]\geq \lambda_n
\end{align*}
holds for every $n\in\mathbb{N}$, and show that $\lim_{n\to\infty}\lambda_n=\lambda$. We omit the details as refer to the proof of \cite[Theorem 1.2]{movingcat}. Altogether, we have shown that
\begin{align*}
\lim_{t\to\infty}\frac{1}{t}\log\mathbb{E}^{(Z,\tilde{\alpha})}_{(0,1)}\left[\exp\left(\int_0^tV(Z(s),\tilde{\alpha}(s))\,\d s\right)\right]=\lambda,
\end{align*}
where, according to the Rayleigh-Ritz formula, $\lambda$ is given by
\begin{align*}
\lambda=&\sup_{f\in\ell^2(\mathbb{Z}^d\times\{0,1\}),\|f\|_2=1}\left<(\tilde{L}+V)f,f\right>.
\end{align*}
Let us calculate the inner product. We have
\begin{align*}
\left<V f, f\right>=-s_0\sum_{x\in\mathbb{Z}^d}f(x,0)^2-s_1\sum_{x\in\mathbb{Z}^d}f(x,1)^2+\gamma f(0,1)^2
\end{align*}
and
\begin{align*}
\left<\tilde{L}f,f\right>=&\sum_{i\in{\{0,1\}}}\sum_{x\in\mathbb{Z}^d}\left((i\kappa+\rho)\sum_{y\sim x}(f(y,i)-f(x,i))f(x,i)+\sqrt{s_0s_1}(f(x,1-i)-f(x,i))f(x,i)\right)
\\=&\sum_{i\in{\{0,1\}}}\sum_{j=1}^d(i\kappa+\rho)\sum_{x\in\mathbb{Z}^d}(f(x+e_j,i)-f(x,i))f(x,i)+(f(x,i)-f(x+e_j,i))f(x+e_j,i)
\\-&\sum_{x\in\mathbb{Z}^d}\sqrt{s_0s_1}(f(x,1)-f(x,0))^2
\\=&-\frac{1}{2}\sum_{i\in{\{0,1\}}}(i\kappa+\rho)\sum_{x,y\in\mathbb{Z}^d, x\sim y}(f(x,i)-f(y,i))^2-\sum_{x\in\mathbb{Z}^d}\sqrt{s_0s_1}(f(x,1)-f(x,0))^2,
\end{align*}
where the factor $\frac{1}{2}$ comes from summing over ordered pairs $(x,y)$. Now, recall from \eqref{dar} that
\begin{align*}
\lim_{t\to\infty}\frac{1}{t}\log\left<U(t)\right>=\sqrt{s_0s_1}+\lambda
\end{align*}
to conclude. 
\end{proof}
\section{Proof of Theorem \ref{annasy3}}
In this chapter, we give a proof for Theorem \ref{annasy3} and consider a dynamic random environment given by a field of independent random walks with equal jump rate $2d\rho$ starting from a Poisson cloud on $\mathbb{Z}^d$ with intensity $\nu$. More precisely, we define the potential $\xi$ to be
\begin{align*}
\xi(x,t)=\sum_{y\in\mathbb{Z}^d}\sum_{j=1}^{N_y}\delta_x(Y_j^y(t)),
\end{align*}
where $N_y$ is a Poisson random variable with intensity $\nu>0$ for each $y\in\mathbb{Z}^d$ and $\{Y_j^y: y\in\mathbb{Z}^d, j=1,\cdots,N_y, Y_j^y(0)=y\}$ is the collection of random walks with jump rate $2d\rho>0$. 

Our first lemma, which provides a more convenient representation of $\left<U(t)\right>$, is an adaptation of \cite[Proposition 2.1]{intermittency} to our setting for switching random walks:
\begin{lemma}
For all $t\geq 0$ and all $\gamma\in[-\infty,\infty)$,
\begin{align}\label{anntraps}
\left<U(t)\right>=\mathbb{E}_{(0,1)}^{(X,\alpha)}\left[\exp\left(\nu\gamma\int_0^t\alpha(s)v_{(X,\alpha)}(X(s),s)\,\emph{d}s\right)\right],
\end{align}
where $v_{(X,\alpha)}(y,t):\mathbb{Z}^d\times[0,\infty)\to\mathbb{R}$ is the solution of
\begin{align}\label{v_(x,a)}
\left\{\begin{array}{lllr}\frac{\emph{d}}{\emph{d} t}v_{(X,\alpha)}(y,t)&=&\rho\Delta v_{(X,\alpha)}(y,t)+\gamma\delta_{(X(t),\alpha(t))}(y,1)v_{(X,\alpha)}(y,t), &t>0,
\\[10pt]v_{(X,\alpha)}(y,0)&=&1\end{array}\right.
\end{align}
conditioned on a fixed realization of $(X,\alpha)$. 
\end{lemma}
\begin{proof}
The proof is similar to the proof of \cite[Proposition 2.1]{intermittency}, but with the additional component $\alpha$. Write $\mathbb{E}_{\nu}$ for the expectation of a Poisson random variable with intensity $\nu$. As in \cite{intermittency}, we integrate out the Poisson system $\xi$ to obtain
\begin{align*}
\left<U(t)\right>&=\left<\mathbb{E}_{(0,1)}^{(X,\alpha)}\left[\exp\left(\gamma\int_0^t\sum_k\delta_{(Y^k(s),1)}(X(s),\alpha(s))\,\d s\right) \right]\right>\nonumber
\\&=\left<\mathbb{E}_{(0,1)}^{(X,\alpha)}\left[\prod_k\exp\left(\gamma\int_0^t\delta_{(Y^k(s),1)}(X(s),\alpha(s))\,\d s\right) \right]\right>\nonumber
\\&=\mathbb{E}_{(0,1)}^{(X,\alpha)}\prod_{y\in\mathbb{Z}^d}\mathbb{E}_{\nu}\mathbb{E}_y^Y\left[\exp\left(\gamma\int_0^t\delta_{(Y(s
),1)}(X(s),\alpha(s))\,\dt d s\right)\right]
\\&=\mathbb{E}_{(0,1)}^{(X,\alpha)}\prod_{y\in\mathbb{Z}^d}\mathbb{E}_{\nu}\left[v_{(X,\alpha)}(y,t)\right].
\end{align*}
Then, taking the expectation with respect to a Poisson random variable yields
\begin{align}\label{anntrapsproof}
\left<U(t)\right>&=\mathbb{E}_{(0,1)}^{(X,\alpha)}\left[\prod_{y\in\mathbb{Z}^d}\sum_n\frac{(\nu v_{(X,\alpha)}(y,t))^n}{n!}e^{-\nu}\right]\nonumber
\\&=\mathbb{E}_{(0,1)}^{(X,\alpha)}\left[\prod_{y\in\mathbb{Z}^d}\exp\left(-\nu(1-v_{(X,\alpha)}(y,t))\right)\right]\nonumber
\\&=\mathbb{E}_{(0,1)}^{(X,\alpha)}\left[\exp\left(-\nu\sum_{y\in\mathbb{Z}^d}w_{(X,\alpha)}(y,t)\right)\right]
\end{align}
for $w_{(X,\alpha)}:=1-v_{(X,\alpha)}$. Note, that for $\gamma<0$, the quantity $v_{(X,\alpha)}(y,t)$ represents the survival probability of $Y$ with start in $y$ up to time $t$, where the (fixed) trajectory of $X$ is seen as a trap, which tries to capture $Y$ with rate $\gamma$ whenever it crosses the latter and if $\alpha$ takes the value $1$ in that moment. For $\gamma>0$, $v_{(X,\alpha)}(y,t)$ represents the number of particles build out of one single particle starting in $0$ which moves around with total jump rate $2d\kappa$ and branches into two, whenever it meets the random walk $X$ and if $\alpha$ equals to $1$ at this time. Next, we see that
\begin{align*}
\sum_{y\in\mathbb{Z}^d}\frac{\dt d}{\dt d t}w_{(X,\alpha)}(y,t)&=-\sum_{y\in\mathbb{Z}^d}\frac{\dt d}{\dt d t}v_{(X,\alpha)}(y,t)
\\&=-\sum_{y\in\mathbb{Z}^d}\left(\rho\Delta v_{(X,\alpha)}(y,t)+\gamma\delta_{(X(t),\alpha(t))}(y,1)v_{(X,\alpha)}(y,t)\right)
\\&=-\gamma\alpha(t)v_{(X,\alpha)}(X(t),t),
\end{align*}
which together with the initial condition $\sum_{y\in\mathbb{Z}^d}w_{(X,\alpha)}(0)=0$ yields
\begin{align*}
\sum_{y\in\mathbb{Z}^d}w_{(X,\alpha)}(y,t)=-\gamma\int_0^t\alpha(s)v_{(X,\alpha)}(X(s),s)\,\dt ds.
\end{align*}
Combined with \eqref{anntrapsproof}, this proofs the claim.
\end{proof}
Before continuing our investigations regarding to the asymptotics of \eqref{anntraps}, we will first consider the case $\kappa=0$, i.\,e.\,, the case of an immobile particle $X$ staying in $0$ the whole time. This idea is highly inspired by \cite{intermittency} and \cite{drewitz} and will be extended here to the case of switching random walks. For $\kappa=0$ the equation \eqref{v_(x,a)} with $X\equiv 0$ reduces to
\begin{align}\label{v_a}
\left\{\begin{array}{llll}
\frac{\d}{\d t}v_{(0,\alpha)}(y,t)&=&\rho\Delta v_{(0,\alpha)}(y,t)+\gamma\alpha(t)\delta_{0}(y)v_{(0,\alpha)}(y,t), &y\in\mathbb{Z}^d, t>0,
\\[10pt]v_{(0,\alpha)}(y,0)&=&1, &y\in\mathbb{Z}^d,\end{array}\right.
\end{align}
such that the annealed survival probability becomes
\begin{align}\label{ann0}
\left<U(t)\right>=\mathbb{E}_{1}^{\alpha}\left[\exp\left(\nu\gamma\int_0^t\alpha(s)v_{(0,\alpha)}(0,s)\,\d s\right)\right].
\end{align}
As we will see later, the following two propositions will help us with the general case $\kappa\geq 0$:
\begin{prop}\label{propk0}
Let $\gamma\in[-\infty,0)$ and $\kappa=0$. Then, as $t\to\infty$,
\begin{align}\label{asytraps0}
\log\left<U(t)\right>=\left\{\begin{array}{ll}\displaystyle-4\nu\sqrt{\frac{\rho s_0}{(s_0+s_1)\pi}}\sqrt{t}(1+o(1)), &d=1,\\[13pt]\displaystyle-4\nu\frac{\rho\pi s_0}{s_0+s_1}\frac{t}{\log\left(t\right)}(1+o(1)), &d=2,\\[13pt]\displaystyle-\tilde{\lambda}_{d,\gamma} t(1+o(1)), &d\geq 3,\end{array}\right.
\end{align}
with
\begin{align*}
\tilde{\lambda}_{d,-\infty}=\inf_{a\in[0,1]}\left\{2d\nu\rho G_d(0)^{-1}a+I(a)\right\},
\end{align*}
where $G_d(0)$ is the Green's function of a simple symmetric random walk in $0$, and
\begin{align*}
\tilde{\lambda}_{d,\gamma}=\inf_{a\in[0,1]}\left\{\frac{2d\nu\rho}{\frac{2d\rho}{|\gamma|}+G_d(0)}a+I(a)\right\}
\end{align*}
for $\gamma\in(-\infty,0)$. 
\end{prop}
\begin{proof}
We start with the case of soft traps $\gamma\in(-\infty,0)$. Recall the representation \eqref{ann0} of the annealed survival probability as well as the solution $v_{(0,\alpha)}(y,t)$ of \eqref{v_a}, which is the survival probability of $Y$ up to time $t$, if we interpret $0$ as a trap which tries to kill $Y$ with rate $|\gamma|$ at time $s$ if $(Y(s),\alpha(s))=(0,1)$. We observe that
\begin{align*}
\int_0^t\alpha(s)\delta_0(Y(s))\d s=\int_0^t\delta_0(Y(s))\d L_s(1)=\int_0^{L_t(1)}\delta_0(Y(S(s))\d s=\int_0^{L_t(1)}\delta_0(\tilde{Y}(s))\d s
\end{align*}
for $S(s):=\inf\{r\geq 0: L_r(1)>s\}$ the right-continuous inverse of $L_\cdot(1)$ and $\tilde{Y}(s):=Y(S(s))$, where we used the time change formula with respect to increasing processes (Note that $\tilde{Y}$ is not Markovian on its own, but Markovian given $\alpha$). Thus, 
\begin{align*}
v_{(0,\alpha)}(0,t)=\mathbb{E}_0^Y\left[\exp\left(\gamma\int_0^t\alpha(s)\delta_0(Y(s))\d s\right)\right]&=\mathbb{E}_0^{Y}\left[\exp\left(\gamma\int_0^{L_t(1)}\delta_0(\tilde{Y}(s))\d s\right)\right]
\end{align*}
For a fixed realization of $\alpha$, the time-change $S(s)$ is deterministic (given $\alpha$) and hence preserves the distributional properties of $Y$. In other words $\tilde{Y}=Y$ in distribution, given $\alpha$. This justifies writing
\begin{align*}
v_{(0,\alpha)}(0,t)=\mathbb{E}_0^Y\left[\exp\left(\gamma\int_0^{L_t(1)}\delta_0(Y(s))\d s\right)\right]=:v(0,L_t(1)),
\end{align*}
where $v$ is independent of $\alpha$ and solves the differential equation
\begin{align}\label{defv}
\left\{\begin{array}{llll}\frac{\d}{\d t}v(y,t)&=&\rho\Delta v(y,t)+\gamma\delta_0 v(y,t), & y\in\mathbb{Z}^d, t>0,
\\[8pt]v(y,0)&=&1, & y\in\mathbb{Z}^d\end{array}\right..
\end{align}
In \cite{drewitz} it has been shown that $v(0,t)$ has the asymptotics
\begin{align}\label{asyv}
v(0,t)=\left\{\begin{array}{ll}\displaystyle\frac{1}{|\gamma|}\sqrt{\frac{\rho}{\pi}}\frac{1}{\sqrt{t}}(1+o(1)), &d=1,\\[12pt]\displaystyle\frac{4\pi\rho}{|\gamma|} \frac{1}{\log(t)}(1+o(1)), &d=2,\\[12pt]\displaystyle\frac{2d\rho}{2d\rho-\gamma G_d(0)}(1+o(1)),&d\geq 3,\end{array}\right.
\end{align}
as $t\to\infty$. 
Therefore
\begin{align*}
\sum_yw_{(0,\alpha)}(y,t)&=-\gamma\int_0^t\alpha(s)v_{(0,\alpha)}(0,s)\d s=-\gamma\int_0^t\alpha(s)v(0,L_s(1))\d s=-\gamma\int_0^{L_t(1)}v(0,u)\d u,
\end{align*}
as $\d L_s(1)=\alpha(s)\d s$. This yields
\begin{align*}
\left<U(t)\right>=\mathbb{E}_1^{\alpha}\left[\exp\left(\nu\gamma \int_0^{L_t(1)}v(0,u)\d u\right)\right].
\end{align*}
where we have the asymptotics
\begin{align*}
\nu\gamma\int_0^{L_t(1)}v(0,s)\d s=\left\{\begin{array}{ll}\displaystyle-4\nu\sqrt{\frac{\rho}{\pi}}\sqrt{L_t(1)}(1+o(1)), &d=1,\\[12pt]\displaystyle-4\nu\pi\rho\frac{L_t(1)}{\log(L_t(1))}(1+o(1)), &d=2,\\[12pt]\displaystyle\nu\gamma\frac{2d\rho}{2d\rho-\gamma G_d(0)}L_t(1)(1+o(1)),&d\geq 3,\end{array}\right.
\end{align*}
for $t\to\infty$. Recall the large deviation principle for the normalized local times $\left(\frac{1}{t}L_t(1)\right)_{t\geq 0}$ of $\alpha$ in state $1$ from Theorem \ref{LDP} with rate function  $I$ given by \eqref{I}, which has a unique zero at $\frac{s_0}{s_0+s_1}$. We can now apply Lemma \ref{varadhanversion} to the functions $f(t)=\sqrt{t}$ and $f(t)=t/\log(t)$ in dimension $1$ and $2$ respectively, to obtain the asymptotics stated in \eqref{asytraps0}. In dimensions $d\geq 3$, Varadhan's Lemma tells us that
\begin{align*}
-\lim_{t\to\infty}\frac{1}{t}\log\mathbb{E}_1^{\alpha}\left[\exp\left(\frac{2d\nu\rho}{\frac{2d\rho}{\gamma}-G_d(0)}L_t(1)\right)\right]=\inf_{a\in[0,1]}\left\{I(a)-\frac{2d\nu\rho}{\frac{2d\rho}{\gamma}-G_d(0)}a\right\}=:\tilde{\lambda}_{d,\gamma}.
\end{align*}
This proofs the proposition for the case $\gamma\in(-\infty,0)$. We proceed with the case of \emph{hard traps}, i.\,e.\,, $\gamma=-\infty$, where the random walk $X$ is immediately killed after crossing one of the traps, if $\alpha$ takes the value $1$ at this times. Then, $v_{(0,\alpha)}$ corresponds to the probability $\phi_{(0,\alpha)}$, that, given $\alpha$, the random walk $Y$ has not hit $0$ up to time $t$ at the same time when $\alpha$ was $0$, i.e.
\begin{align*}
\phi_{(0,\alpha)}(y,t)=\mathbb{P}_y^Y(\forall s\leq t: (Y(s),\alpha(s))\neq (0,1)).
\end{align*}
An analogous argument as above shows that
\begin{align*}
\phi_{(0,\alpha)}(y,t)=\mathbb{P}_y^Y(\forall s\leq L_t(1):Y(s)\neq 0)=\phi(y,L_t(1))
\end{align*}
given $\alpha$, where $\phi(y,t):=\mathbb{P}_y^Y(\forall s\leq t: Y(s)\neq 0)$. Set 
\begin{align*}
\psi_{\alpha}(y,t):=1-\phi_{(0,\alpha)}(y,t)=\mathbb{P}_y^Y(\exists s\in[0,t]:Y(s)=0,\alpha(s)=1),
\end{align*}
which plays the role of $w_{(0,\alpha)}$ and so
\begin{align*}
\left<U(t)\right>=\mathbb{E}_1^{\alpha}\left[\exp\left(-\nu\sum_{y\in\mathbb{Z}^d}\psi_{\alpha}(y,t)\right)\right].
\end{align*}
Further, let $\psi(y,t):=1-\phi(y,t)$ denote the probability, that $Y$ with start in $y$ hits $0$ at least once up to time $t$, regardless of values of $\alpha$, which solves the differential equation
\begin{align*}
\left\{\begin{array}{llll}
\frac{\d}{\d t}\psi(y,t)&=&\rho\Delta\psi(y,t), &t>0, y\neq 0,\\[10pt]\psi(0,t)&=&1, &t>0,\\[10pt]\psi(y,0)&=&0, &y\neq 0.
\end{array}\right.
\end{align*}
Due to the relation of $\psi$ and $\phi$ we have $\psi_\alpha(y,t)=\psi(y,L_t(1)$ given $\alpha$, and therefore, due to the chain rule,
\begin{align*}
\frac{\d}{\d t}\psi(y,L_t(1))=\frac{\d}{\d t}L_t(1)\cdot\frac{\d}{\d t}\psi(y,\cdot)(L_t(1))=\alpha(t)\cdot \Delta\psi(y,L_t(1)).
\end{align*}
Hence,
\begin{align*}
\frac{\d}{\d t}\sum_{y\in\mathbb{Z}^d}\psi_{\alpha}(y,t)&=\frac{\d}{\d t}\sum_{y\in\mathbb{Z}^d}\psi(y,L_t(1))=\sum_{y\in\mathbb{Z}^d}\alpha(t)\Delta\psi(y,L_t(1))-\alpha(t)\Delta\psi(0,t)
\\&=-\alpha(t)\rho\sum_{y\sim 0}(\psi(y,L_t(1))-\psi(0,L_t(1))=-2d\rho\alpha(t)(\psi(e_1,L_t(1))-1)
\end{align*}
with $e_1=(1,0,\cdots,0)^T$ the first unit vector, where we used the symmetry of the random walk as well as the fact that $\sum_y\Delta \psi_{\alpha}(y,t)=0$. Thus,
\begin{align}\label{intlt}
\sum_{y\in\mathbb{Z}^d}\psi_{\alpha}(y,t)=\int_{0}^t2d\rho\phi(e_1,L_s(1))\alpha(s)\d s=\int_0^{L_t(1)}2d\rho\phi(e_1,s)\d s
\end{align}
where we substituted $L_s(1)$ in the last step. Now, the quantity $\phi(e_1,t)$ is known (see e.g.\,\cite{lawler}) to have the asymptotics
\begin{align*}
\phi(e_1,t)=\left\{\begin{array}{ll}\sqrt{\frac{1}{\pi\rho t}}(1+o(1)), &d=1,\\[12pt]\frac{\pi}{\log(t)}(1+o(1)), &d=2,\\[12pt]G_d(0)^{-1}(1+o(1)),&d\geq 3,\end{array}\right.
\end{align*}
as $t\to\infty$, where $G_d$ is the Green's function of a $d$-dimensional symmetric random walk with generator $\Delta$. Thus,
\begin{align}\label{psia}
\sum_{y\in\mathbb{Z}^d}\psi_{\alpha}(y,t)=\left\{\begin{array}{ll}4\sqrt{\frac{\rho}{\pi}}\sqrt{L_t(1)}(1+o(1)), &d=1,\\\\4\pi\rho \frac{L_t(1)}{\log(L_t(1))}(1+o(1)), &d=2,\\\\2d\rho G_d(0)^{-1}L_t(1)(1+o(1)),&d\geq 3,\end{array}\right.
\end{align}
as $t\to\infty$, given $\alpha$. As the function $g\colon [0,1]\to\mathbb{R}$, $x\mapsto -\nu\rho G_d(0)^{-1}x$ is continuous and bounded, we can apply Varadhan's lemma to deduce the limit
\begin{align*}
-\lim_{t\to\infty}\frac{1}{t}\log\mathbb{E}_1^{\alpha}\left[\exp\left(-\nu\frac{2d\rho}{G_d(0)}L_t(1)\right)\right]=\inf_{a\in[0,1]}\left\{I(a)+\frac{2d\nu\rho}{G_d(0)}a\right\}=:\tilde{\lambda}_{d,\infty}.
\end{align*}
This establishes the asymptotics \eqref{asytraps0} for $d\geq 3$. For $d=1$, we apply lemma \ref{varadhanversion} with $f(t)=\sqrt{t}$ to obtain
\begin{align*}
\lim_{t\to\infty}\frac{1}{\sqrt{t}}\log\mathbb{E}_1^{\alpha}\left[\exp\left(-\nu\sqrt{\frac{8\rho}{\pi}}\sqrt{L_t(1)}\right)\right]=-4\nu\sqrt{\frac{\rho}{\pi}}\sqrt{\frac{s_0}{s_0+s_1}}.
\end{align*}
The case $d=2$ is similar with $f(t)=t/\log(t)$.
\end{proof}
Note that in the first two dimensions the survival probability decays sub-exponentially and does not depend on $\gamma$ and in higher dimensions $d\geq 3$ the asymptotics for $\gamma=\infty$ are consistent with those of the case $\gamma<\infty$, as $\lim_{\gamma\to\infty}\tilde{\lambda}_{d,\gamma}=\tilde{\lambda}_{d,\infty}$.
\\\\The next proposition deals with the case $\gamma>0$ of catalysts, still under the assumption $\kappa=0$.
\begin{prop}\label{propk02}
Let $\gamma\in(0,\infty)$ and $\kappa=0$. Then for all dimensions $d\geq 1$ the annealed number of individuals grows with double-exponential rate
\begin{align*}
\lim_{t\to\infty}\frac{1}{t}\log\log\left<U(t)\right>=\sup_{f\in\ell^2(\mathbb{Z}^d),\|f\|_2=1}\left(\gamma f(0)^2-\frac{1}{2}\sum_{x,y\in\mathbb{Z}^d, x\sim y}\rho(f(x)-f(y))^2\right).
\end{align*}
\end{prop}
\begin{proof}
Recall the representation \eqref{anntraps} as well as the solution $v_{(0,\alpha)}$ to \eqref{v_a} for a fixed realization of $\alpha$, which can also be written as
\begin{align*}
v_{(0,\alpha)}(0,t)=\mathbb{E}_0^Y\left[\exp\left(\gamma\int_0^t\alpha(s)\delta_0(Y(s))\d s\right)\right]
\end{align*}
in point $0$. Let
\begin{align*}
v(0,t):=\mathbb{E}_0^Y\left[\exp\left(\gamma\int_0^t\delta_0(Y(s))\d s\right)\right]
\end{align*}
which plays the analogous role of the function $v$ in the proof of Proposition \ref{propk0}, but now for positive $\gamma$. In \cite{movingcat} it has been shown that
\begin{align}\label{grenzwert}
\lim_{t\to\infty}\frac{1}{t}\log v(0,t)=\mu
\end{align}
where
\begin{align*}
\mu:=\sup_{f\in\ell^2(\mathbb{Z}^d),\|f\|_2=1}\left(\gamma f(0)^2-\frac{1}{2}\sum_{x,y\in\mathbb{Z}^d, x\sim y}\rho(f(x)-f(y))^2\right)
\end{align*}
is the largest eigenvalue of the self-adjoint operator $\mathcal{H}:=\rho\Delta+\gamma\delta_0$ on $\ell^2(\mathbb{Z}^d)$. Furthermore, it is known from \cite{intermittency} that $\mu$ is always positive in dimensions $d=1,2$ and
\begin{align*}
\mu\left\{\begin{array}{ll}=\displaystyle 0, &\displaystyle0<\frac{\gamma}{\rho}\leq \frac{1}{G_d(0)},
\\\displaystyle >0, &\displaystyle\frac{\gamma}{\rho}> \frac{1}{G_d(0)}
\end{array}\right.
\end{align*}
in dimensions $d\geq 3$, where $G_d$ is the Green's function of a simple symmetric random walk with jump rate $2d$. Next, we aim to compare $v(0,t)$ to $v_{(0,\alpha)}(0,t)$ in the same manner as in the proof of Proposition \ref{propk0}, which yields
\begin{align}\label{v_0zeit+}
v_{(0,\alpha)}(0,t)=v(0,L_t(1)).
\end{align}
By the spectral theorem for $\mathcal{H}$ and the simplicity of $\mu$, 
\begin{align*}
v(0,t)=\e^{\mu t}(1+o(1)), \qquad t\to\infty,
\end{align*}
uniformly in $t$ large. Hence,
\begin{align*}
\int_0^{L_t(1)}v(0,s)\,\d s=\int_0^{L_t(1)}\e^{\mu s}(1+o(1))\,\d s, \qquad t\to\infty
\end{align*}
and therefore
\begin{align*}
\left<U(t)\right>&=\mathbb{E}_{1}^{\alpha}\left[\exp\left(\nu\gamma\int_0^{L_t(1)}v(0,s)\,\d s\right)\right]\\&=\mathbb{E}_1^{\alpha}\left[\exp\left(\nu\gamma\int_0^{L_t(1)}\e^{\mu s}(1+o(1))\,\d s\right)\right]
\\&=\mathbb{E}_1^{\alpha}\left[\exp\left(\frac{\nu\gamma}{\mu}\e^{\mu L_t(1)}(1+o(1))\right)\right]
\end{align*}
as $t\to\infty$, where the error does not depend on the fluctuations of $\alpha$, as seen in a similar manner than before. Now, on one hand, we have the upper bound
\begin{align*}
\left<U(t)\right>\leq \exp\left(\frac{\nu\gamma}{\mu}\e^{\mu t}(1+o(1))\right), \quad t\to\infty,
\end{align*}
and on the other hand,
\begin{align}\label{theother}
\mathbb{E}_1^{\alpha}\left[\exp\left(\frac{\nu\gamma}{\mu}\e^{\mu L_t(1)}(1+o(1))\right)\right]\geq& \exp\left(\frac{\nu\gamma}{\mu}\e^{\mu t}(1+o(1))\right)\mathbb{P}_1^{\alpha}(L_t(1)=t),
\end{align}
as $t\to\infty$. Now, recall the large deviation principle for for the normalized local times with rate function $I$ on scale $t$, which asserts that
\begin{align}\label{loglog}
\lim_{t\to\infty}\frac{1}{t}\log\mathbb{P}_1^{\alpha}(L_t(1)=t)=-I(1).
\end{align}
Hence, combining \eqref{loglog} with \eqref{theother},
\begin{align*}
\lim_{t\to\infty}\frac{1}{t}\log\log\mathbb{E}_1^{\alpha}\left[\exp\left(\frac{\nu\gamma}{\mu}\e^{\mu L_t(1)}(1+o(1))\right)\right]\geq& \lim_{t\to\infty}\frac{1}{t}\log\log\left(\exp\left(\frac{\nu\gamma}{\mu}\e^{\mu t}\right)\mathbb{P}_1^{\alpha}(L_t(1)=t)\right)
\\=&\lim_{t\to\infty}\frac{1}{t}\log\left(\frac{\nu\gamma}{\mu}\e^{\mu t}+\log\mathbb{P}_1^{\alpha}(L_t(1)=t)\right)
\\=&\lim_{t\to\infty}\frac{1}{t}\log\left(\frac{\nu\gamma}{\mu}\e^{\mu t}\right)=\mu.
\end{align*}
\end{proof}
The next ingredient for the proof of our main result is the so called \emph{Pascal principle} which asserts that if we average over the environment, than the best the random walk $X$ can do in order to maximize its mass is to stay still in the starting point, which brings us back to the case $\kappa=0$. In the setting of a simple random walk without any switching component, this has been proven in \cite{drewitz} and  \cite{intermittency} for $\gamma<0$ and $\gamma>0$ respectively. Therefore, the question arises if the Pascal principle can still provide an upper bound also in our case of a switching random walk, or if there is a better joint strategy of the random walk together with the dormancy component $\alpha$. The next lemma ensures that the Pascal principle still proves best also in our case.
\begin{lemma}\label{pascalpoisson} Recall $v_{(X,\alpha)}$ as the solution to \eqref{v_(x,a)} and the solution $v_{(0,\alpha)}$ of \eqref{v_a} for any fixed realization of $(X,\alpha)$. Then, for all $\gamma\in[-\infty,0)$, $y\in\mathbb{Z}^d$ and $t\geq 0$,
\begin{align*}
\int_0^t\alpha(s)v_{(X,\alpha)}(X(s),s)\,\emph{d}s\geq \int_0^t\alpha(s)v_{(0,\alpha)}(X(s),s)\,\emph{d}s
\end{align*}
and for all $\gamma\in(0,\infty)$, $y\in\mathbb{Z}^d$ and $t\geq 0$,
\begin{align*}
v_{(X,\alpha)}(y,t)\leq v_{(0,\alpha)}(0,t).
\end{align*}
\end{lemma}
\begin{proof}
First, let $\gamma\in[-\infty,0)$ and recall that $v_{(0,\alpha)}(y,t)=v(y,L_t(1))$. Now, \cite[Proposition 2.1]{drewitz} asserts that for any piecewise constant function $\hat{X}:[0,t]\to\mathbb{Z}^d$ with a finite number of discontinuities,
\begin{align*}
\int_0^tv_{\hat{X}}(\hat{X}(s),s)\,\dt ds\geq \int_0^tv(0,s)\,\dt ds
\end{align*}
where $v_{\hat{X}}$ is the solution to
\begin{align}\label{defhatX}
\left\{\begin{array}{llll}\frac{\d}{\d t}v_{\hat{X}}(y,t)&=&\rho\Delta v_{\hat{X}}(y,t)+\gamma\delta_{\hat{X}(t)}(y)v_{\hat{X}}(y,t),&y\in\mathbb{Z}^d, t>0,\\[10pt]v_{\hat{X}}(y,0)&=&1, &y\in\mathbb{Z}^d.
\end{array}\right.
\end{align}
In a similar way as before, we see that if $v_{(\hat{X},\alpha)}$ is the solution to
\begin{align*}
\left\{\begin{array}{llll}\frac{\d}{\d t}v_{(\hat{X},\alpha)}(y,t)&=&\rho\Delta v_{(\hat{X},\alpha)}(y,t)+\gamma\alpha(t)\delta_{\hat{X}(t)}(t)v_{(\hat{X},\alpha)}(y,t),&y\in\mathbb{Z}^d, t>0,
\\[10pt]v_{(\hat{X},\alpha)}(y,0)&=&1, &y\in\mathbb{Z}^d,\end{array}\right.
\end{align*}
then $v_{(\hat{X},\alpha)}(y,t)=v_{\hat{X}}(y,L_t(1))$, as $\hat{X}$ is independent of $\alpha$, and hence,
\begin{align*}
\int_0^t\alpha(s)v_{(\hat{X},\alpha)}(\hat{X}(s),s)\,\d s&=\int_0^t\alpha(s)v_{\hat{X}}(\hat{X}(s),L_s(1))\,\d s=\int_0^{L_t(1)}v_{\hat{X}}(\hat{X}(s),s)\,\d s
\\&\geq \int_0^{L_t(1)}v(0,s)\d s=\int_0^t\alpha(s)v(0,L_s(1))\d s\\&=\int_0^t\alpha(s)v_{(0,\alpha)}(0,s)\d s.
\end{align*}
Finally, conditioned on $\alpha$, the random walk $X$ is a piecewise constant function with a finite number of discontinuities, such that we can replace $\hat{X}$ with $X$.

Next, let $\gamma>0$. The proof works along the same lines as in \cite[Proposition 2.2]{intermittency}. From this proof we already know that if $p_{\rho}(x,y)$ denotes the probability density function of a random walk with generator $\rho\Delta$ and start in $0$, then
\begin{align*}
\max_{x\in\mathbb{Z}^d}p_{\rho}(x,y)=p_{\rho}(0,t)
\end{align*}
for all $t\geq 0$. Further, in \cite[Proposition 2.2]{intermittency} has been shown that $h^{*n}\to 0, n\to\infty$, uniformly on compact intervals, where $h^{*n}$ denotes the n-fold convolution of the function
\begin{align*}
h(t):=\gamma p_{\rho}(0,t).
\end{align*}
Since
\begin{align*}
v_{(X,\alpha)}(X(t),t)=&1+\gamma\int_0^tp_{\rho}(X(t)-X(s),t-s)\alpha(s)v_{(X,\alpha)}(X(s),s)\,\d s
\\\leq &1+\gamma\int_0^tp_{\rho}(0,t-s)\alpha(s)v_{(X,\alpha)}(X(s),s)\,\d s
\end{align*}
and
\begin{align*}
v_{(0,\alpha)}(0,t)=&1+\gamma\int_0^tp_{\rho}(0,t-s)\alpha(s)v_{(0,\alpha)}(0,s)\,\d s,
\end{align*}
we have
\begin{align}\label{h1}
v_{(X,\alpha)}(X(\cdot),\cdot)\leq 1+h^*(\alpha v_{(X,\alpha)}(X(\cdot),\cdot))
\end{align}
as well as
\begin{align}\label{h2}
v_{(0,\alpha)}(0,\cdot)= 1+h^*(\alpha v_{(0,\alpha)}(0,\cdot)).
\end{align}
Hence,
\begin{align*}
v_{(0,\alpha)}(0,\cdot)-v_{(X,\alpha)}(X(\cdot),\cdot)\geq h^{*n}(\alpha(v_{(0,\alpha)}(0,\cdot)-v_{(X,\alpha)}(X(\cdot),\cdot))
\end{align*}
by iteration and substraction of \eqref{h1} and  \eqref{h2}. As $h^{*n}\to 0$ for  $n\to\infty$, we can deduce
\begin{align*}
v_{(X,\alpha)}(X(\cdot),\cdot)\leq v_{(0,\alpha)}(0,\cdot)
\end{align*}
for all realizations of $(X,\alpha)$. Altogether, we obtain
\begin{align*}
v_{(X,\alpha)}(y,t)=&1+\gamma\int_0^tp_{\rho}(y-X(s),t-s)\alpha(s)v_{(X,\alpha)}(X(s),s)\,\d s
\\\leq &1+\gamma\int_0^tp_{\rho}(0,t-s)\alpha(s)v_{(X,\alpha)}(X(s),s)\,\d s
\\\leq &1+\gamma\int_0^tp_{\rho}(0,t-s)\alpha(s)v_{(0,\alpha)}(0,s)\d s=v_{(0,\alpha)}(0,t),
\end{align*}
as desired. 
\end{proof}
We are now ready to prove Theorem 1.3.
\begin{proof}[Proof of Theorem \ref{annasy3}]
The proof can be summarized as follows: We will first show that the case $\kappa=0$ considered in the Proposition \ref{propk0} and \ref{propk02} provides a lower bound for the general case $\kappa\geq 0$ in dimensions $d=1,2$. This gives us together with the upper bound asserted in Lemma \ref{pascalpoisson} the desired asymptotics. 

Let us start with the case $\gamma<0$ of traps and show that the asymptotics of $\left<U(t)\right>$ are lower-bounded by \eqref{asytraps0} in dimensions $d\in\{1,2\}$. We adapt the notations from the proof of \cite[Lemma 2.1]{drewitz}: Let $E_t$ be the event that none of the traps starts in a ball $B_{R_t}$ of radius $R_t$ around $0$, where we choose $R_t$ to be $\small\frac{t}{\ln(t)}$ for $d=1$ resp.\,$\ln(t)$ for $d=2$. This event has the probability $\e^{-\nu(R_t+1)^d}$. Further, let $G_t$ be the event that $X$ with start in $0$ stays in $B_{R_t}$ up to time $t$. Analogously, we define $\tilde{G}_t$ to be the event that a simple symmetric random walk $\tilde{X}$ without switching and with jump rate $2d\kappa$ stays in $B_{R_t}$ up to time $t$. Then, using time-change again,
\begin{align*}
\mathbb{P}(G_t)&=\mathbb{P}_{(0,1)}^{(X,\alpha)}(X(s)\in B_{R_t}\forall s\leq t)=\mathbb{P}_1^{\alpha}\mathbb{P}_0^{\tilde{X}}(\tilde{X}(s)\in B_{R_t}\forall s\leq L_t(1))
\\&\geq \mathbb{P}_1^{\alpha}\mathbb{P}_0^{\tilde{X}}(\tilde{X}(s)\in B_{R_t}\forall s\leq t)=\mathbb{P}(\tilde{G}(t))\geq \exp\left(\ln(\beta)\frac{t}{R_t^2}\right),
\end{align*}
where the last inequality is known from \cite{drewitz} for some $\beta>0$. Moreover, let $F_t$ be the event that each trap which starts outside $B_{R_t}$ only intersects $B_{R_t}$ during time periods where $\alpha$ takes the value $0$, i.e.\,, where $X$ is dormant. Then, we can lower-bound
\begin{align*}
\left<U(t)\right>\geq \mathbb{P}(E_t)\mathbb{P}(F_t)\mathbb{P}(G_t).
\end{align*}
Note that the event $F_t$ differs from the analogous event appearing in the proof of \cite[Lemma 2.1]{drewitz}, which was the event that each trap which starts outside $B_{R_t}$ never enters $B_{R_t}$ up to time $t$. Here, making use of the protection provided by the dormancy mechanism, we can relax this condition and only require that the traps stay outside of $B_{R_t}$ whenever $X$ is active. In order to compute $\mathbb{P}(F_t)$ we first look at $\tilde{F}_t$ which shall denote the event that no trap $Y$ which starts in $y\neq 0$ ever hits $0$ at time points $s$ with $\alpha(s)=1$. The probability of this event is nothing but the survival probability of $X$ in case $\kappa=0$, which has been studied in Proposition \ref{propk0}. Thus, $\mathbb{P}(\tilde{F}_t)$ is asymptotically equal to \eqref{asytraps0}. Comparing $\mathbb{P}(F_t)$ to $\mathbb{P}(\tilde{F}_t)$ exactly in the same way as in \cite[Lemma 2.1]{drewitz}, we find that these are asymptotically equal. Finally, we compare the decay rates of all probabilities $\mathbb{P}(E_t),\mathbb{P}(F_t),\mathbb{P}(G_t)$ for $t\to\infty$ to conclude that the annealed survival probability is asymptotically lower bounded by $\mathbb{P}(F_t)$ and hence \eqref{asytraps0}. 

Note, that we did not prove a lower bound in dimension $d=3$, so that we can only deduce the existence of a constant $\lambda_{d,\gamma,\rho,\kappa,s_0,s_1}\geq \tilde{\lambda}_{d,\gamma}$ which may depend on all the parameters.

We continue with the case $\gamma>0$. The upper bound again follows from the Pascal principle stated in Lemma \ref{pascalpoisson}. For the lower bound, we force the random walk $X$ to stay in the starting point $0$ up to time $t$ and use the fact that, for a simple symmetric random walk $\tilde{X}$ without switching and with jump rate $2d\kappa$, 
\begin{align*}
\mathbb{P}_0^{\tilde{X}}(\tilde{X}(s)=0\forall s\in[0,t])=\e^{-2d\kappa t}.
\end{align*}
Moreover, recall from Proposition \ref{propk02} that
\begin{align*}
\left<U_0(t)\right>=\exp\left(\frac{\nu\gamma}{\mu}\e^{\mu t}(1+o(1))\right), \quad t\to\infty,
\end{align*}
where $U_0(t)$ shall denote the number of individuals up to time $t$ in the case $\kappa=0$. Hence, 
\begin{align*}
\lim_{t\to\infty}\frac{1}{t}\log\log\left<U(t)\right>&\geq
\lim_{t\to\infty}\frac{1}{t}\log\log\left(\exp\left(\frac{\nu\gamma}{\mu}\e^{\mu t}(1+o(1))\right)\cdot\e^{-2d\kappa t}\cdot\mathbb{P}_1^{\alpha}(\alpha(s)=1\forall s\in[0,t])\right)
\\&=\lim_{t\to\infty}\frac{1}{t}\log\left(\frac{\nu\gamma}{\mu}\e^{\mu t}(1+o(1))-2d\kappa t-I(1)t\right)
\\&=\mu,
\end{align*}
where we used the large deviation principle for the local times of $\alpha$ with rate function $I$. 
\end{proof}
\subsection*{Acknowledgement}
The author would like to thank Professor Wolfgang König for his invaluable support and many helpful suggestions. 

\subsection*{Competing Interests}
The author has no competing interests to declare that are relevant to the content of this article.

\end{document}